\documentclass[12pt,psamsfonts]{amsart}

\newtheorem{anyprop}{Anyprop}[section]

\newtheorem{theorem}[anyprop]{Theorem}
\newtheorem{lemma}[anyprop]{Lemma}
\newtheorem{proposition}[anyprop]{Proposition}
\newtheorem{corollary}[anyprop]{Corollary}

\theoremstyle{definition}

\newtheorem{definition}[anyprop]{Definition}
\newtheorem{example}[anyprop]{Example}

\newtheorem{question}[anyprop]{Question}

\newtheorem{remark}[anyprop]{Remark}

\theoremstyle{remark}

\numberwithin{equation}{section}

\newcommand{\NN}{\mathbb{N}}
\newcommand{\ZZ}{\mathbb{Z}}

\newcommand{\PP}{\mathbb{P}}
\renewcommand{\AA}{\mathbb{A}}

\newcommand  {\shE}     {\mathcal{E}}
\newcommand  {\shF}     {\mathcal{F}}

\newcommand  {\shM}     {\mathcal{M}}

\newcommand  {\shS}     {\mathcal{S}}
\newcommand  {\shT}     {\mathcal{T}}

\newcommand  {\shQ}     {\mathcal{Q}}

\newcommand  {\fom}     {\mathfrak{m}}

\newcommand  {\Char}    {\operatorname{char}}

\newcommand  {\dual}    {\vee}

\newcommand  {\Hom}     {\operatorname{Hom}}

\renewcommand  {\ker }  {\operatorname{ker}}

\newcommand  {\lra}     {\longrightarrow}

\renewcommand{\O}       {\mathcal{O}}

\newcommand  {\Proj}    {\operatorname{Proj}}

\newcommand  {\ra}      {\rightarrow}

\newcommand  {\rk}    {\operatorname{rk}}

\newcommand  {\sign}    {\operatorname{sign}}

\newcommand  {\Spec}    {\operatorname{Spec}}

\newcommand  {\Syz}     {\operatorname{Syz}}

\newcommand{\comdots}{ , \ldots , }
\newcommand{\komdots}{ , \ldots , }
\newcommand{\plusdots}{ + \ldots + }
\newcommand{\oplusdots}{ \oplus \ldots \oplus }

\newcommand{\geqdots}{ \geq \ldots \geq }
\newcommand{\leqdots}{ \leq \ldots \leq }
\newcommand{\capdots}{ \cap \ldots \cap }
\newcommand{\cupdots}{ \cup \ldots \cup }
\newcommand{\wedgedots}{ \wedge \ldots \wedge }

\newcommand{\uplusdots}{ \uplus \ldots \uplus }

\newcommand{\tensor}{\otimes}

\newcommand{\indi}{i}

\newcommand{\indj}{j}

\newcommand{\indk}{k}

\newcommand{\indsec}{{j}}

\newcommand{\numiii}{\renewcommand{\labelenumi}{(\roman{enumi})}}

\newcommand{\fug}{{g}}
\newcommand{\fuh}{{h}}

\newcommand{\point}{{P}}

\newcommand{\new}{\newcommand}

\newcommand{\buntor}{{\shE}} %toricbundle
\newcommand{\subbuntor}{{\shF}} %toricbundle
\newcommand{\fibspec}{{E}} %Special fiber W
\newcommand{\subfibspec}{{F}} %subspace U
\newcommand{\vecspec}{{w}}

\new{\tup}{{\nu}}
\new{\tupmon}{{\sigma}}
\new{\indgen}{{k}}
\new{\varex}{{z}}
\new{\sect}{{z}}
\new{\eles}{{u}}
\new{\coox}{{x}}
\new{\foli}{{\psi}}
\new{\folitor}{{\psi}}
\new{\seco}{{;}}
\new{\tupsub}{{\nu}}

\newcommand{\dreidreieckmatrix}[4]
{\left(
\begin{array}{ccc}
#1 & \ldots & #2 \cr
\ldots & \ldots & \ldots \cr
#3 & \ldots & #4
\end{array}
\right)}

\newcommand{\dreisechseckmatrix}[8]
{\left(
\begin{array}{cccccc}
#1 & \ldots & #2 & \ldots & #3 & #4 \cr : & : & : & : &: &: \cr #5 &
\ldots & #6 & \ldots & #7 & #8
\end{array}
\right)}

\def\mydate{\number\day\space\ifcase\month \or January\or February\or March\or April\or May\or
June\or July\or August\or September\or October\or November\or
December\fi \space\number\year}

\usepackage{amscd}
\usepackage{amssymb}

\newcommand{\fieldk}{K}

\newcommand{\hcf}{{\rm hcf}}

\newcommand{\indpart}{{\lambda}}
\newcommand{\numpart}{{\ell}}

\usepackage{ifthen}

\newboolean{paper}
\setboolean{paper}{false}
\begin{document}

\title[Looking out for stable syzygy bundles]
{Looking out for stable syzygy bundles}

\author[Holger Brenner]{Holger Brenner}

\maketitle

With an appendix by Georg Hein: Semistability of the general syzygy
bundle.

%\thanks{}

%\date{}

% at present the "communicated by" line appears only in ERA and PROC
%\commby{}

%\dedicatory{Preliminary version,  \mydate}

\begin{abstract}
We study (slope-)stability properties of syzygy bundles on a projective space $\PP^N$
given by ideal generators of a homogeneous primary ideal.
In particular we give a combinatorial criterion for a monomial ideal
to have a semistable syzygy bundle.
Restriction theorems for semistable bundles yield the same stability results
on the generic complete intersection curve.
From this we deduce a numerical formula for the tight closure of an ideal
generated by monomials or by generic homogeneous elements
in a generic two-dimensional complete intersection ring.
\end{abstract}

\noindent
Mathematical Subject Classification (2000):
13A35; 13D02; 14D20; 14H60; 14J60;

%===========================================================
\section*{Introduction}

Suppose that $f_1 \komdots f_{d+1} \in \fieldk[U,V]=R$ are $d+1$
generic homogeneous polynomials of degree $d$ in the two-dimensional
polynomial ring over a field $\fieldk$. Since the dimension of the
space of forms of degree $d$ is $d+1$, it follows that these
generically chosen elements form a basis, and therefore we get the
ideal inclusion $R_{ \geq d}  \subseteq (f_1 \komdots f_{d+1})$ and
hence also
$$  R_{ \geq d+1}  \subseteq (f_1 \komdots f_{d+1}) \,\,\,\, \,\,\, (*)   \, .$$
This last statement is by no means true for other two-dimensional
standard-graded domains such as $R=\fieldk[X,Y,Z]/(G)$. One aim of
this paper is to show that $(*)$ is true for such a two-dimensional
hypersurface ring (for $G$ generic of sufficiently high degree), if
we replace the ideal  $(f_1 \komdots f_{d+1})$ on the right hand
side by its tight closure $(f_1 \komdots f_{d+1})^*$ (Corollary
\ref{tightgenerate}). This means that $d+1$ generic forms of degree
$d$ are ``tight generators'' for $R_{\geq d+1}$.

The theory of tight closure has been developed by M. Hochster and C. Huneke since 1986
and plays a central role in commutative algebra (\cite{hochsterhunekebriancon}, \cite{hunekeapplication},
\cite{hunekeparameter}).
It assigns to every ideal $I$ in a Noetherian ring containing a field an ideal $I^* \supseteq I$,
which is called the tight closure of $I$.
For a domain over a field of positive characteristic $p$ it is defined with the help of the Frobenius
endomorphism, by
$$I^* := \{ f\in R :\, \exists c \neq 0 \mbox{ such that } c f^q \in I^{[q]}  \mbox{ for all } q=p^{e} \} \, .$$

The tight closure of an ideal in a regular ring is the ideal itself, and
it is a typical feature of this theory that we may generalize
a statement about an ideal in a regular ring
to a non-regular ring if we replace the ideal by its tight closure.
The tight closure version of the Theorem of Brian\c{c}on-Skoda is an important instance for this principle,
and our stated result fits well into this picture.

There are three main ingredients for the above mentioned result and for similar results in this paper:

1) The geometric interpretation of tight closure and slope criteria.

2) Restriction theorems for stable vector bundles.

3) Criteria for stable syzygy bundles on a projective space.

\medskip
We explain in this introduction these three items and their
interplay and we give a summary on the content of this paper.

\medskip
1) Geometric interpretation of tight closure and slope criteria

We will use the geometric approach to the theory of tight closure in
terms of vector bundles which we have developed in
\cite{brennertightproj}, \cite{brennerslope} and
\cite{brennerhabil}. The starting point is the cohomological
characterization of tight closure due to Hochster saying that $f \in
(f_1 \komdots f_n)^*$ holds for an $\fom$-primary ideal $(f_1
\komdots f_n)$ in an excellent normal local domain $(R,\fom)$ of
dimension $d$ over a field of positive characteristic if and only if
$H_\fom^d (A) \neq 0$, where $A=R[T_1 \komdots T_n]/(f_1T_1
\plusdots f_nT_n+f)$ is the forcing algebra for these
data.\footnote{A remark about the characteristic: the theory arising
in characteristic $0$ from this cohomological characterization is
called solid closure; see \cite{hochstersolid}. However solid
closure has in dimension two all the good properties which we expect
for a tight closure type theory and we will take it in this paper as
the technical definition of tight closure and denote it henceforth
with $\star $.}

If $R$ is a normal two-dimensional standard-graded domain over an
algebraically closed field $\fieldk$ and the data $f_1 \komdots f_n$
and $f$ are $R_+$-primary and homogeneous, then this cohomological
characterization takes a simple form in terms of the locally free
sheaf of syzygies $\Syz(f_1 \komdots f_n)$ on the smooth projective
curve $C =  \Proj \, R$. This syzygy bundle is given by the short
exact sequence
$$ 0 \lra \Syz(f_1 \comdots f_n)(m)
\lra \bigoplus_{i=1}^n \O_C(m-\deg (f_i)) \stackrel{f_i}{\lra}
\O_C(m) \lra 0 \, .$$ In this situation $f \in (f_1 \komdots
f_n)^\star $ holds if and only if the affine-linear bundle
corresponding to the cohomology class $c=\delta(f) \in H^1(C,
\Syz(f_1 \komdots f_n)(m))$ (where $m =\deg (f)$) is not an affine
scheme.

This geometric approach allows us to apply the elaborated methods of
algebraic geometry to problems coming from tight closure. In
\cite{brennerslope} we studied the ampleness and bigness properties
of the dual of the syzygy bundle $\Syz(f_1 \komdots f_n)(m)$ in
dependence of $m$ and obtained both inclusion and exclusion criteria
for tight closure in terms of the minimal and the maximal slope of
it. These criteria together yield under the condition that the
syzygy bundle is semistable (we shall recall the definitions in
section \ref{restrictionsection}) the numerical characterization
that
$$ (f_1 \komdots f_n)^\star = (f_1 \komdots f_n) + R_{\geq \frac{\deg(f_1) \plusdots \deg(f_n)}{n-1}} \, $$
holds in characteristic $0$ (see \cite[Theorem 8.1]{brennerslope}
and Remark \ref{positive} below for results in positive
characteristic).

In order to apply this numerical formula to the computation of tight
closure one has to establish the semistability property of a given
syzygy bundle on the projective curve $C = \Proj\, R$. This is a
difficult matter in general, even if the rank of the bundle is $2$
(and the number of ideal generators is $3$). One result of
\cite{brennercomputationtight} is that the syzygy bundle
$\Syz(X^d,Y^d,Z^d)$ is semistable on normal domains
$R=\fieldk[X,Y,Z]/(G)$ for $\deg (G) \geq 3d-1$ and therefore $
(X^d,Y^d,Z^d)^\star =(X^d,Y^d,Z^d) + R_{\geq 3d/2}$.

Another and more general way to obtain semistable syzygy bundles on
curves is to work on the projective plane (or a projective space or
other varieties in which the curve lives) and to apply restriction
theorems.

\medskip
2) Restriction theorems for stable vector bundles.

There exist beautiful theorems due to Mehta and Ramanathan, Flenner,
Bogomolov and Langer (see \cite{mehtaramanathanrestriction},
\cite{flennerrestriction}, \cite{bogomolovstability},
\cite{bogomolovstablesurface}, \cite{huybrechtslehn},
\cite{langersemistable}) saying that the restriction of a
(semi)stable bundle on a smooth projective variety $X$ to a general
complete intersection curve of sufficiently high degree is again
(semi)stable. We will present these theorems in section
\ref{restrictionsection}.

We shall use mainly the easiest instance of this type of results,
the restriction of stable bundles on the projective plane $\PP^2$ to
a generic curve $C \subset \PP^2$. Homogeneous elements $f_1
\komdots f_n \in \fieldk[X, Y,Z]$ which are primary to the
irrelevant ideal $(X,Y,Z)$ define a locally free syzygy bundle
$\Syz(f_1 \komdots f_n)$ on $\PP^2$ and its restriction to a
projective curve $C = V_+(G)$ gives the bundle which is crucial for
the computation of the solid closure $(f_1 \komdots f_n)^\star$ in
$R=\fieldk[X,Y,Z]/(G)$. So if we know that the syzygy bundle
$\Syz(f_1 \komdots f_n)$ is semistable on $\PP^2$, the restriction
theorems yield at once that the same is true for $\Syz(f_1 \komdots
f_n)|_C$ for a general curve $C$ of sufficiently high degree. This
gives us then the generic answer for $(f_1 \komdots f_n)^\star$ in a
two-dimensional hypersurface ring. The result of Flenner gives a
bound for the degree of the curve and the result of Bogomolov shows
moreover that the restriction is in fact semistable for every smooth
curve fulfilling a stronger degree condition.

So we are led to look out for stable syzygy bundles on the projective plane or more generally on a projective space.
Note that the restriction theorems give us the possibility to argue on a regular polynomial ring to obtain
results on tight closure, though ``tight closure does nothing''(Hochster) on regular rings!

\medskip
3) Criteria for stable syzygy bundles on a projective space.

Our main problem is now: suppose that homogeneous elements $f_1
\komdots f_n \in \fieldk[X_0 \komdots X_N]$ are given. When is the
syzygy bundle $\Syz(f_1 \komdots f_n)$ on $\PP^N$ semistable? The
sections \ref{numericalsection} - \ref{generalsection} are concerned
with this question.

There exist surprisingly few results on stability properties of
syzygy bundles. Flenner shows in \cite[Corollary
2.2]{flennerrestriction} (also proved by Ballico in \cite[Corollary
6.5]{ballicorestriction}) that the syzygy bundle for all monomials
of fixed degree is semistable.

In section \ref{numericalsection} we shall discuss necessary conditions for a syzygy bundle to be semistable.
We get results by comparing the slope of $\Syz(f_1 \komdots f_n)$ with the slopes of the natural subsheaves
$\Syz (f_i, i \in J)$ for subfamilies $J \subset \{1 \komdots n\}$.
This gives at once the necessary degree condition
$d_n \leq \frac{d_1 \plusdots d_{n-1}}{n-2}$ for semistability, where $d_n$ is the largest degree
(Proposition \ref{numcondition}).

The stability of the syzygy bundle implies conditions for the
existence of global sections of the bundle and of its dual. These
observations provide easily a characterization of semistability for
bundles of rank $2$ and $3$, which correspond to $n=3$ and $4$ ideal
generators. We can take advantage of the fact that there exist only
few line bundles on a projective space, contrary to the situation on
projective curves (section \ref{sectionlowrank}).

In section \ref{linesection} we study the restriction of a syzygy bundle on $\PP^N$
to generic lines $\PP^1 \subset \PP^N$.
If these restrictions are a direct sum of line bundles of the same degree, then the bundle itself is semistable.
Since this property is fulfilled for $d+1$ generic forms of degree $d$, their syzygy bundle is semistable.
Hence we may derive the result mentioned at the beginning of the introduction (Proposition \ref{restrictionline},
Corollary \ref{tightgenerate}).

A torsion free subsheaf $\shT \subseteq \Syz(f_i, i \in I)$ of rank $r$
yields an invertible subsheaf
$(\bigwedge^r \shT)^{\dual \dual} \subseteq \bigwedge^r (\Syz ( f_i, i \in I))$.
Therefore we deal in section
\ref{wedgesection}
with exterior powers of syzygy bundles and describe them as a kernel
of a suitable mapping.

In sections \ref{monomialsection} and \ref{examples} we settle the
case of the syzygy bundle of a monomial ideal using results of
A. Klyachko on toric bundles (\cite{klyachkoselecta}, \cite{klyachkospectralproblems},
\cite{klyachkoequivariant}).
The main result is that $\Syz(X^{\tupmon_i}, i \in I)$ is semistable
if for every subset
$J \subseteq I$ the corresponding subsheaf $\Syz(X^{\tupmon_i}, i
\in J) \subseteq \Syz(X^{\tupmon_i}, i \in I)$ does not contradict
the semistability. This provides an easy combinatorial test for
semistability in the monomial case.

Finally, section \ref{generalsection} addresses the case of ideals
which are generated by generic forms $f_i$ of degrees $d_i$
fulfilling the necessary numerical conditions from section
\ref{numericalsection}. From a Theorem of Bohnhorst-Spindler we
deduce that the syzygy bundle of $n$ parameters in an
$n$-dimensional polynomial ring is semistable (Corollary
\ref{bohnhorstcorollary}) and a Theorem of Hein asserts that this is
also true for the syzygy bundle of $n$ generic forms of degree $d$
under the condition that $n \leq d(N+1)$. This theorem is proven by
Hein in an appendix to this paper.

I thank the referee and Almar Kaid for critical remarks.

\section{Stable bundles and restriction theorems}
\label{restrictionsection}

We recall the definition of semistability on a smooth projective
curve $C$ over an algebraically closed field $\fieldk$. Let $\shS$
denote a locally free sheaf on $C$ of rank $r$. The degree of $\shS$
is defined as the degree of the corresponding invertible sheaf $
\det \shS= \bigwedge^r \shS$. The number $\mu (\shS) = \deg ( \shS)
/\rk (\shS)$ is called the slope of the vector bundle. A locally
free sheaf $\shS$ is called semistable, if for every locally free
subsheaf $\shT \subset \shS$ the inequality $\mu( \shT) \leq \mu
(\shS)$ holds (and stable if $<$ holds). This notion is due to
Mumford \cite{mumfordprojective} and plays a crucial role in the
construction of moduli spaces for vector bundles on curves and
beyond.

On a higher dimensional smooth projective variety it is convenient to develop
these notions more generally for torsion-free coherent sheaves
$\shS$ in dependence of a fixed very ample invertible sheaf. We will
work here only with the notion of $\mu$-stability (or
Mumford-Takemoto stability), not with Gieseker stability. We take
\cite{okonekschneiderspindler} as our main reference and we deal
only with coherent torsion-free sheaves on a projective space
$\PP^N$.
The determinantal bundle of such a sheaf is defined by the bidual
$\det \shS = (\bigwedge^r \shS) ^{\dual \dual} $, which is an invertible sheaf,
and the degree of $\shS$ is defined by $\deg (\bigwedge^r \shS) ^{\dual \dual}$.
Since $\shS$ is locally free outside a closed subset of codimension $\geq 2$, there exist projective lines $\PP^1 \subset \PP^N$
such that the restriction is locally free,
hence $\shS|{\PP^1} \cong \O(a_1) \oplusdots \O(a_r)$
and this gives another way to define the degree, as $\sum a_i$.
The slope of $\shS$ is defined as before be dividing through the rank.

\begin{definition}
Let $\shS$ denote a torsion-free coherent sheaf on a projective space $\PP^N$.
Then $\shS$ is called semistable if
for every coherent subsheaf $\shT \subseteq \shS$ the inequality
$\mu(\shT) \leq \mu(\shS)$ holds.
\end{definition}

These subsheaves are of course again torsion-free. It is enough to check this property for those subsheaves
which have a torsion-free quotient (see \cite[Theorem 1.2.2]{okonekschneiderspindler}).

The restriction of a semistable torsion-free sheaf to a curve is in general not semistable anymore.
We will use the following restriction theorems which we cite here for the convenience of the reader.
We only state the results for a bundle on a projective space and for the restriction to a complete
intersection curve.

\begin{theorem}
\label{mehtaramanathanrestriction} {\rm(}Mehta-Ramanathan, see
\cite[Theorem 6.1]{mehtaramanathanrestriction}, \cite[Theorem
7.2.8]{huybrechtslehn}{\rm)} Let $\fieldk$ denote an algebraically
closed field of any characteristic and let $\shS$ denote a
semistable torsion-free coherent sheaf on $\PP^N$. Then there exists
a number $k_0$ such that for $N-1$ general elements $D_1 \komdots
D_{N-1} \in |\O(k)|$, $k \geq k_0$, the restriction $\shS|C$ is
again semistable on the smooth complete intersection curve $C=D_1
\capdots D_{N-1}$.
\end{theorem}

This Theorem of Mehta-Ramanathan says nothing about the bound $k_0$. This is provided by the Theorem of Flenner,
but only in characteristic zero.

\begin{theorem}
\label{flennerrestriction} {\rm(}Flenner, see
\cite{flennerrestriction},\cite[Theorem 7.1.1]{huybrechtslehn}{\rm)}
Let $\fieldk$ denote an algebraically closed field of characteristic
$0$. Let $\shS$ denote a torsion-free coherent semistable sheaf of
rank $r$ on the projective space $\PP^N$. Then for $k$ fulfilling
the condition that
$$ \frac{ \binom{k+N}{N} -(N-1)k -1}{k} > \max \{ \frac{r^2-1}{4}, {1} \}$$
and for $N-1$ general elements $D_1 \komdots D_{N-1} \in |\O(k)|$
the restriction $\shS|C$ is again semistable on the smooth complete intersection curve
$C=D_1 \capdots D_{N-1}$.
\end{theorem}

\begin{remark}
\label{rankdegreeremark}
The degree bound in the Theorem of Flenner reduces for $N=2$ to the condition that
$k > \frac{r^2-3}{2}$ and $k \geq 2$.
So this means $k \geq 2$ for $r=2$, $k \geq 4$ for $r=3$, $k \geq 7$ for $r=4$.
\end{remark}

In the surface case the Theorem of Bogomolov
gives even a result about the restriction
to every smooth curve, not only to a general curve.

\begin{theorem}
\label{bogomolovrestriction} {\rm(}Bogomolov, see
\cite{bogomolovstability}, \cite{bogomolovstablesurface},
\cite[Theorem 7.3.5]{huybrechtslehn}{\rm)} Let $\fieldk$ denote an
algebraically closed field of characteristic $0$. Let $\shS$ denote
a stable locally free sheaf of rank $r$ on the projective plane
$\PP^2$ with Chern classes $c_1$ and $c_2$. Let $ \triangle (\shS) =
2rc_2- (r-1) c_1^2$ and set $R= \binom{r}{ \lfloor r/2 \rfloor}
\binom{r-2}{\lfloor r/2 \rfloor -1}$. Then the restriction $\shS|C$
is again stable for every smooth curve $C \subset \PP^2$ of degree
$k$ with $2k > \frac{R}{r} \triangle (S) + 1 $.
\end{theorem}

\begin{remark}
If $\shS= \Syz(f_1 \komdots f_n)$ for polynomials of degree $d_i$,
then $c_1 (\shS) = - \sum d_i$
and $c_2(\shS) = \frac{(\sum d_i)^2 - \sum d_i^2}{2}$.
Therefore the discriminant is in this case
$$\triangle (\shS)= (n-1)( ( \sum d_i)^2- \sum d_i^2) -(n-2)(\sum d_i)^2= (\sum d_i)^2- (n-1) \sum d_i^2 \, .$$
If all the degrees are constant, then $\triangle (\shS)=nd^2$ and Bogomolovs result yields the degree condition
$2k >3d^2 + 1$ for $n=3$,
$2k >4d^2 + 1$ for $n=4$ and
$2k >60d^2 + 1$ for $n=5$.
\end{remark}

\begin{example}
Look at the syzygy bundle $\Syz(X^2,Y^2,Z^2)$ on $\PP^2$. It is easy
to see that this bundle is stable (see Corollary
\ref{threeregularsemistable} below) and the Theorem of Flenner
(\ref{flennerrestriction}) tells us that the restriction is
semistable for a generic curve $C$ of degree $\deg C \geq 2$, $C =
\Proj R$, $R=\fieldk[X,Y,Z]/(G)$. The bound in the Theorem of
Bogomolov tells us that the restriction to every smooth curve of
degree $\geq 7$ is semistable. Due to \cite[Proposition
6.2]{brennercomputationtight} this is even true for degree $ \geq
5$.

For $\deg G=3$ the semistability depends on the curve equation $G$,
and so does the question whether $XYZ \in (X^2,Y^2,Z^2)^\star$ holds
in $R=\fieldk[X,Y,Z]/(G)$ or not. For the Fermat cubic $G=
X^3+Y^3+Z^3$, the semistability property is easy to establish, since
this curve equation gives at once the relation $(X,Y,Z)$ for
$(X^2,Y^2,Z^2)$ (of total degree $3$), which yields a short exact
sequence
$$ 0 \lra \O \lra \Syz(X^2,Y^2,Z^2)(3) \lra \O \lra 0 \, . $$
This shows that the syzygy bundle is semistable (and strongly
semistable, but not stable). It follows that $XYZ \in
(X^2,Y^2,Z^2)^\star$ holds in $\fieldk[X,Y,Z]/(X^3+Y^3+Z^3)$ in any
characteristic, which was first shown by a quite complicated
computation of A. Singh; see \cite{singhcomputation}.
\end{example}

\begin{remark}
\label{positive}
We comment on the situation in positive
characteristic. The best restriction theorem for semistable bundles
in positive characteristic is due to A. Langer
\cite{langersemistable} and gives a Bogomolov-type restriction
theorem. However, the numerical formula for tight closure mentioned
in the introduction needs the assumption that the syzygy bundle is
strongly semistable, meaning that every Frobenius pull-back of it is
semistable. It was shown in \cite{brennerstronglysemistable} that a
Bogomolov-type restriction theorem for strong semistability does not
hold. It is open whether there exists a Flenner-type restriction
theorem for strong semistability.

However, we may derive a slightly weaker result for prime
characteristic $p \gg 0$ from the results in characteristic zero. If
we know in the relative situation, that is over $\Spec \ZZ$, that a
syzygy bundle $\Syz(f_1 \komdots f_n)$ is semistable in
characteristic zero, then every twist of it of positive degree is
ample, and therefore this property holds also in positive
characteristic for almost all prime characteristics $p$. From this
it follows for $p$ large enough that for $\deg(f)  < \frac{\deg(f_1)
\plusdots \deg(f_n)}{n-1}$ the element $f$ belongs to $(f_1 \komdots
f_n)^*$ only if it belongs to the ideal itself, and for $\deg(f)  >
\frac{\deg(f_1) \plusdots \deg(f_n)}{n-1}$ (not $\geq$ as in the
formula) the element belongs to the Frobenius closure of $(f_1
\komdots f_n)$, hence also to the tight closure. In particular, if
the degree bound is not a natural number, then we get the same
statement as in characteristic zero.
\end{remark}

\section{Numerical conditions for semistability on $\PP^N$}
\label{numericalsection}

Let $f_i$, $ i \in J$, denote homogeneous polynomials $\neq 0$ in
$\fieldk[X_0 \comdots X_N]$, where $\fieldk$ is an algebraically
closed field. Their syzygy sheaf $\Syz(f_i, i \in J)$ is locally
free on $D_+(f_i, i \in J)= \bigcup _{i \in J} D_+(f_i)
\subseteq
\PP^N$ and has rank $r=|J|-1$. We compute first the degree of
$\Syz(f_i, i \in J)(m)$, which is by definition the degree of the
invertible sheaf $\det (\Syz(f_i, i \in J)(m))$, where $\det
(\Syz(f_i, i \in J)(m)) =(\bigwedge^r (\Syz(f_i, i \in J)(m))^{\dual
\dual}$.

\begin{lemma}
\label{detidentify}
Let $f_i \in \fieldk[X_0 \komdots X_N]$, $i \in
J$, denote homogeneous polynomials $\neq 0$ of degree $d_i$,
$|J|=r+1$, $r \geq 1$. Then the following hold.

\renewcommand{\labelenumi}{(\roman{enumi})}

\begin{enumerate}

\item
Suppose that the $f_i$ do not have a common factor. Then
$$ \det (\Syz(f_i, i \in J) (m)) \cong \O(rm - \sum_{i \in J} d_i) $$
and $\deg (\Syz(f_i, i \in J) (m)) = rm - \sum_{i \in J} d_i$.

\item
Suppose that the $f_i$ do have a highest common factor $f$ of
degree $d$. Then
$$ \det (\Syz(f_i, i \in J) (m)) \cong \O(rm+d - \sum_{i \in J} d_i) $$
and $\deg (\Syz(f_i, i \in J) (m)) = r m + d - \sum_{i \in J} d_i$.

\ifthenelse{\boolean{paper}}{}{
\item
An identification {\rm(}outside a subset of codimension $\geq
2${\rm)} for the determinantal bundle $\det (\Syz( f_i, i \in
J))$ is for any fixed $k$ on $D_+(g_k)$, $g_k= f_k/f$, given by
$$ s_1 \wedgedots s_{r} \longmapsto \sign(k,J)
\frac{f}{f_k} \det \big( (s_{ji})_{j=1 \komdots r,\,\, i \in J -\{k
\}} \big) \, . $$ {\rm(}Here $J =\{ 1 \komdots r+1 \}$ is supposed
to be ordered and $\sign(k,J)$ is $1$, if $k$ is an even element in
$J$, and $-1$ if it is an odd element.{\rm)}}
\end{enumerate}
\end{lemma}

\begin{proof} (i). Suppose first that the $f_i$, $i \in J$, do not have a
common factor, so that their zero locus $V_+(f_i, i \in J)$ has
codimension $\geq 2$. Thus we have the short exact (presenting)
sequence
$$ 0 \lra \Syz(f_i, i \in J)(m) \lra \bigoplus_{i \in J} \O(m-d_i) \lra \O(m) \lra 0 \, $$
on $D_+(f_i, i \in J)$.
We restrict this sequence to a projective line $\PP^1 \subset D_+(f_i,\, i \in J)$
and get
$$ \deg (\Syz(f_i, i \in J) (m)) = (r+1) m - \sum_{i \in J}  d_i -m  \, .$$

(ii) Now suppose that the $f_i$ do have a highest common factor,
and write $f_i=f g_i$ such that the $g_i$, $i \in J$, do not have a
common factor. Then we have an isomorphism of sheaves
$$ \Syz(f_i, i \in J)(m) \cong \Syz(g_i, i \in J)(m - d)  $$
by considering a syzygy $(s_i, i \in J)$ for $g_i$ of total degree
$m-d$ as a syzygy for $f_i=fg_i$ of total degree $m$. Therefore we
get
\begin{eqnarray*}
\deg (\Syz(f_i, i \in J)(m))
&=& \deg (\Syz(g_i, i \in J)(m-d)) \cr
&=& r( m - d) - \sum _{i \in J} (d_i -d) \cr
\ifthenelse{\boolean{paper}}{}{&=& rm     -rd - \sum_{i \in J} d_i + (r+1) d
\cr}
&=& rm - \sum_{i \in J} d_i + d  \, .
\end{eqnarray*}
\ifthenelse{\boolean{paper}}{\end{proof}}{
(iii). For a fixed $k \in J$ we consider now the mapping
$$ \bigwedge^r (\Syz(f_i))
\lra \bigwedge^r (\bigoplus_{i \in J, i \neq k} \O(-d_i)) \cong\O(-
\sum_{i \in J, i \neq k} d_i) \stackrel{\sign (k,J)
\frac{f}{f_k}}{\lra } \O(d - \sum_{i \in J} d_i) \, .$$ This mapping
sends the wedge product of $r$ syzygies $s_j, j=1 \komdots r$ of
total degree $0$ to
$$
s_1 \wedgedots s_r \longmapsto  \sign(k,J) \frac{f}{f_k} \det
\big((s_{ji}) _{j=1 \komdots r,\,\, i \in J -\{ k\}} \big) \, .$$
This mapping is well-defined on $\bigcup_{i \in J} D_+(g_i)$, since
for $k \leq r$ ($\check{s_{j,k}}$ means omit this)
\begin{eqnarray*}
& & \frac{\sign(k,J)}{f_k}
\det \dreisechseckmatrix  {s_{1,1}}  {\check{s_{1,k}}} {s_{1,r}} {s_{1,r+1}} {s_{r,1}}  {\check{s_{r,k}}} {s_{r,r}} {s_{r,r+1}}  \cr
&= & \frac{\sign(k,J)}{f_k}
\det \dreisechseckmatrix  {s_{1,1}}  {\check{s_{1,k}}} {s_{1,r}} {\frac{-1}{f_{r+1}} \sum_{i=1}^r f_is_{1,i}}
{s_{r,1}}  {\check{s_{r,k}}} {s_{r,r}} {\frac{-1}{f_{r+1}} \sum_{i=1}^r f_is_{r,i}  }  \cr
&= & \frac{\sign(k,J)}{f_k}
\det \dreisechseckmatrix  {s_{1,1}}  {\check{s_{1,k}}} {s_{1,r}} {\frac{-f_k}{f_{r+1}}s_{1,k}}
{s_{r,1}}  {\check{s_{r,k}}} {s_{r,r}} {\frac{-f_k}{f_{r+1}}s_{r,k}}  \cr
&= & \frac{-\sign(k,J)}{f_{r+1}}
 \det \dreisechseckmatrix  {s_{1,1}}  {\check{s_{1,k}}} {s_{1,r}} {s_{1,k}}
{s_{r,1}}  {\check{s_{r,k}}} {s_{r,r}} {s_{r,k}}  \cr
&= & \frac{-\sign(k,J)}{f_{r+1}}  (-1)^{r-k}
\det \dreidreieckmatrix  {s_{1,1}} {s_{1,r}}  {s_{r,1}} {s_{r,r}}     \, .
\end{eqnarray*}
This fits well together, since the sign is now $(-1) (-1)^k
(-1)^{r-k} = (-1)^{r+1} = \sign(r+1,J)$. This mapping sends on $D_+(f_1 \komdots f_{r+1})$ the wedge product
$$(\frac{1}{f_{1}} ,0 \komdots 0, \frac{-1}{f_{r+1}} )
\wedgedots (  0   \komdots    0 , \frac{1}{f_{r}} , \frac{-1}{f_{r+1}}  )$$
to
$$ \sign(r+1,J)   f/f_1 \cdots f_{r+1}  \neq 0 \, .$$
So this must be an isomorphism, since it is an endomorphism of an
invertible sheaf on $\bigcup_{i \in J} D_+(g_i)$.
\end{proof}}

This gives the following necessary numerical conditions for a sheaf of syzygies to be semistable.

\begin{proposition}
\label{subsetcondition}
Let $f_i \neq 0$, $i \in I$, $|I| \geq 2$,
denote homogeneous elements in the polynomial ring $\fieldk[X_0
\comdots X_N]$ of degrees $d_i$. For every subset $J \subseteq I$
denote by $d_J$ the degree of the highest common factor of the
subfamily
$f_i$, $i \in J$. Suppose that the syzygy sheaf
$\Syz(f_i, i \in I)$ is semistable. Then for every $J \subseteq I$,
$|J| \geq 2$, we have the numerical condition
$$ \frac{d_J- \sum_{i \in J} d_i}{|J|-1}  \leq \frac{d_I- \sum_{i \in I} d_i}{|I|-1}   \; .  $$
If $\Syz(f_i, i \in I)$ is stable, then $<$ holds for $J \subset I$.
\end{proposition}
\begin{proof} Every subset $J \subseteq I$ defines the syzygy subsheaf
$\Syz(f_i, i \in J)$ sitting in
$$ 0 \lra \Syz(f_i, i \in J) \lra  \Syz(f_i, i \in I) \lra \bigoplus_{i \not\in J} \O(-d_i) \, .$$
(The sequence is not exact on the right in general.)
The semistability of $\Syz(f_i, i \in I)$ implies that
$$ \mu (\Syz(f_i, i \in J)) \leq \mu (\Syz(f_i, i \in I)) \; ,$$
and we have computed these slopes in Lemma \ref{detidentify} (ii).
\end{proof}

\begin{remark}
This necessary condition for semistability is in general not sufficient, as Example
\ref{necnotsuff} below shows.
However, if the $f_i$ are monomials, then we will see in Section \ref{monomialsection}
(Corollary \ref{monomstablecrit})
that this condition is also sufficient.
\end{remark}

The condition in Proposition \ref{subsetcondition}
implies the following necessary condition for the degrees of a
semistable syzygy sheaf.

\begin{corollary}
\label{numcondition}
Let $f_1 \komdots f_n \in \fieldk[X_0 \comdots X_N]$ denote homogeneous polynomials $\neq 0$ without common factor
of degrees $1 \leq d_1 \leqdots d_n$. Suppose that their syzygy
sheaf is semistable. Then for every $1 \leq r \leq n-2$ we have the
numerical condition
$$ (n-r-1) (d_1 + \ldots +d_{r+1}) \geq r(d_{r+2} + \ldots +d_n) \; .  $$
For $r=n-2$
this gives the necessary condition
$ d_1  \plusdots d_{n-1} \geq (n-2)d_n $.
On the other hand, this last condition implies the other ones.
\end{corollary}

\begin{proof}
We apply Proposition \ref{subsetcondition} to the subset
$J= \{1 \komdots r+1 \} $ and get the inequality
$$ \frac{- \sum _{i =1}^{r+1} d_i}{r}  \leq \frac{- \sum _{i =1}^n d_i}{n-1}$$
or equivalently that
$$ (n-1) \sum_{i=1}^{r+1} d_i  \geq r \sum_{i =1}^n d_i \, .$$
Subtracting
$r \sum_{i=1}^{r+1} d_i$ gives the result.

For the last statement, suppose that $d_1  \plusdots d_{n-1} \geq
(n-2) d_n $ holds and that we have proved already that $(n-1)
\sum_{i=1}^{r+1} d_i  \geq r \sum_{i=1}^n d_i$ (descending induction
on $r$). We have $(n-1)d_{r+1} \leq (n-1) d_n \leq \sum_{i=1}^n
d_i$. Therefore
$$ (n-1) \sum_{i=1}^{r} d_i
= (n-1) \sum_{i=1}^{r+1}d_i - (n-1)d_{r+1}
\geq r \sum_{i=1}^n d_i - (n-1)d_{r+1}
 \geq (r-1) \sum_{i=1}^n d_i \, .$$
\end{proof}

\section{Syzygy bundles of low rank}
\label{sectionlowrank}

We cover now the case of a syzygy sheaf of rank $2$ and $3$
(corresponding to $3$ or $4$ ideal generators). The following
criteria are known.

\begin{lemma}
\label{ranktwothreecrit}
Suppose that $\shS$ is a coherent
torsion-free sheaf on a projective space $\PP^N$ over a field. Then
the following criteria for semistability hold.
\renewcommand{\labelenumi}{(\roman{enumi})}
\begin{enumerate}

\item
Suppose that $\rk (\shS)=2$.
If $\shS(m)$ has no global sections $\neq 0$ for $m < -\mu(\shS) = - \deg(\shS)/2$, then $\shS$
is semistable.

\item
Suppose that $\rk (\shS)=3$.
Suppose that $\shS(m)$ has no global sections $\neq 0$ for $m < -\mu (\shS ) = -\deg(\shS)/3$
and that $(\shS^\dual)(k)$ has no global sections $\neq 0$ for $k < - \mu (\shS^\dual) = \deg(\shS)/3$.
Then $\shS$ is semistable.
\end{enumerate}
\end{lemma}
\begin{proof}
See \cite[Lemma 1.2.5, Remark 1.2.6]{okonekschneiderspindler}.
\end{proof}

This gives the following corollaries for $3$ ideal generators.

\begin{corollary}
\label{threeregularsemistable} Let $f_1,f_2,f_3 \in \fieldk[X_0
\komdots X_N]$ be a homogeneous regular sequence with degrees $d_1
\leq d_2 \leq d_3$ such that $d_3 \leq d_1 +d_2$. Then the sheaf
$\Syz(f_1,f_2,f_3)$ is semistable on $\PP^N$ {\rm(}and stable for
$<${\rm )}.
\end{corollary}
\begin{proof}
The Koszul complex yields the resolution
$$0 \lra \O(m-d_1-d_2-d_3) \lra  \bigoplus_{i \neq j} \O(m-d_i -d_j) \lra \Syz (m) \lra 0 \, .$$
Since $N \geq 2$, the global sections of $\Syz(m) $ come from the left,
hence $\Gamma(\PP^N, \Syz(m))=0$ for $m < d_1+d_2$.
So $\Syz(m)$ has no non-trivial sections for $m < \frac{d_1+d_2+d_3}{2}  \leq d_1+d_2 $
and the result follows from Lemma \ref{ranktwothreecrit}.
\end{proof}

\begin{corollary}
\label{threeprimarydegree} Let $\fieldk$ denote an algebraically
closed field of characteristic $0$. Let $f_1,f_2,f_3 \in \fieldk[X,
Y,Z]$ be homogeneous primary elements with degrees $d_1 \leq d_2
\leq d_3$ such that $d_3 \leq d_1 +d_2$. Then for a general
hypersurface ring $R=\fieldk[X, Y,Z]/(G)$ for general homogeneous
$G$ of degree $ \geq 2$ we have
$$ (f_1,f_2,f_3)^\star = (f_1,f_2,f_3) + R_{\geq \frac{d_1+d_2+d_3}{2}} \, .$$
The same is true for every $G$ defining a smooth curve under
the degree condition $ \deg (G) \geq d_1d_2+d_1d_3+d_2d_3 - \frac{1}{2}(d_1^2+d_2^2+d_3^2) +1 $.
\end{corollary}
\begin{proof}
This follows from Corollary \ref{threeregularsemistable} and the restriction theorem
of Flenner \ref{flennerrestriction}.
The last statement follows from the restriction theorem of Bogomolov \ref{bogomolovrestriction}, since
$\triangle (\Syz) = 2d_1d_2+2d_1d_3+2d_2d_3 -d_1^2-d_2^2-d_3^2$
and $R/r = 1$.
\end{proof}

\begin{example}
A typical example where we may apply Corollary
\ref{threeprimarydegree} is for the computation of
$(X^d,Y^d,Z^d)^\star$. The generic answer is that
$(X^d,Y^d,Z^d)^\star =(X^d,Y^d,Z^d) + R_{\geq 3d/2}$ and this holds
for $R=\fieldk[X,Y,Z]/(G)$ for general $G$ of degree $\deg G \geq 2$
and for every normal $R$ for $\deg G \geq \frac{3}{2} d^2 +1$.
\end{example}

\begin{example}
Consider the elements $X^3$, $XY^2$ and $ZY^2$ in $\fieldk[X,Y,Z]$.
These polynomials are not $(X,Y,Z)$-primary and their syzygy sheaf
is locally free only outside the points $(0,0,1)$ and $(0,1,0)$. It
fulfills the degree condition in Corollary \ref{threeprimarydegree},
but it is not semistable. The syzygy $(0,Z,-X)$ is a non-trivial
global section of $\Syz(X^3,XY^2,ZY^2)(4)$, but its degree is $2
\cdot 4 -9=-1$ negative.
\end{example}

We consider now the case of $4$ polynomials in three variables.

\begin{corollary}
\label{fourprimaryelements} Let $f_1,f_2,f_3,f_4 \in \fieldk[X,
Y,Z]$ be homogeneous primary elements with ordered degrees $d_1 \leq
d_2 \leq d_3 \leq d_4$. Suppose that $2d_4 \leq d_1 +d_2+d_3$ and
that $\Gamma(\PP^2, \Syz(m))= 0$ for $m <
\frac{d_1+d_2+d_3+d_4}{3}$. Then the syzygy bundle
$\Syz(f_1,f_2,f_3,f_4)$ is semistable.

Furthermore, for $\Char (K)=0$ and for a general hypersurface ring
$R=\fieldk[X, Y,Z]/(G)$ for $G$ of degree $ \geq 4$ we have
$$ (f_1,f_2,f_3,f_4)^\star = (f_1,f_2,f_3,f_4) + R_{\geq \frac{d_1+d_2+d_3+d_4}{3}} \, .$$
The same is true for every $G$ defining a smooth curve and fulfilling
the degree condition $\deg (G) \geq \sum_{i \neq j}d_id_j-\sum_i d_i^2 +1 $.
\end{corollary}
\begin{proof}
We dualize the presenting sequence and get
$$0 \lra \O (-m)  \lra \bigoplus_{i=1}^4 \O( d_i -m) \lra  (\Syz(m))^\dual \lra 0 \, .$$
Since $H^1(\PP^2, \O(-m))=0$ for all $m \in \ZZ$, every global
section of $ (\Syz(m))^\dual$ comes from a global section of
$\bigoplus_{i=1}^4 \O(m- d_i)$. This means that every cosyzygy
$\Syz(m) \ra \O$ must factor through $\bigoplus_{i=1}^4 \O(m- d_i)$.
Therefore for $m
> d_4$ there exists no non-trivial homomorphism $\Syz(m) \ra \O$.
From the assumption it follows that $ \frac{d_1+d_2+d_3+d_4}{3} \geq
d_4$, hence there exists no cosyzygy for $m >
\frac{d_1+d_2+d_3+d_4}{3}$. So both conditions in Lemma
\ref{ranktwothreecrit} for $\shS= \Syz(0)$ hold true and the result
follows.

The statements about solid closure follows from the Theorems of
Flenner \ref{flennerrestriction} and Bogomolov
\ref{bogomolovrestriction}. \end{proof}

\begin{example}
\label{necnotsuff} Consider the four elements $X^{10},
Y^{10},Z^{10}$ and $P=X^9Y+X^9Z+Y^9X+Y^9Z+Z^9X+Z^9Y$. All the syzygy
subbundles $\Syz(f_i,i \in J)$ for subsets $J \subset \{1,2,3,4\}$
do not contradict the semistability. This is clear for $|J|=2$ since
the polynomials are pairwise coprime and for $|J|=3$ since the
numerical degree condition in Corollary \ref{fourprimaryelements} is
fulfilled. However, this syzygy bundle is not semistable. We have
$XYZP \in (X^{10},Y^{10},Z^{10})$ and therefore there exists a
non-trivial syzygy of degree $13$. But the degree of $\Syz(13)$ is
$3 \cdot 13 -4 \cdot 10=-1$.
\end{example}

\section{Semistable restrictions to a generic projective line}
\label{linesection}
Let $\shS$ denote a coherent torsion-free sheaf
$\shS$ on $\PP^N$. The slope of $\shS$ and of a subsheaf $\shT
\subseteq \shS$ can be computed on a generic line $\PP^1 \subset
\PP^N$. Hence if we know that the restriction of $\shS$ to a generic
projective line is semistable, that is of type $\O(a) \oplusdots
\O(a)$, then $\shS$ is semistable (see \cite[Remark after Lemma
2.2.1]{okonekschneiderspindler}). We derive from this observation
the following semistability result for $d+1$ forms of degree $d$ and
we obtain the result mentioned in the introduction that $d+1$
general elements of degree $d$ are ``tight generators'' for the next
degree in a generic two-dimensional complete intersection ring.

\begin{proposition}
\label{restrictionline} Let $f_1 \komdots f_{d+1} \in \fieldk[X_0
\komdots X_N]$ be $d+1$ forms of degree $d$ over an algebraically
closed field $\fieldk$. Suppose that there exists a linear morphism
$\fieldk[X_0 \komdots X_N] \ra \fieldk[U,V]$ such that the images of
these forms are linearly independent in
$\fieldk[U,V]_d$. Then the syzygy bundle $\Syz(f_1, \ldots ,
f_{d+1})$ is semistable on $\PP^N$. This is in particular true for
$d+1$ generic forms of degree $d$.
\end{proposition}

\begin{proof}
Let the linear mapping be given by $X_j \mapsto a_jU+b_jV$. We may
write the images as $\tilde{f}_i = \sum_{k=0}^d  p_{i,k}(a_j,b_j)
U^kV^{d-k}$, where the coefficients $p_{i,k}$ are polynomials in
$a_j,b_j$. Since the determinant of the $(d+1) \times (d+1)$ matrix
of this polynomial entries is $\neq 0$ for a special value of
$(a_j,b_j)$, the determinant is not the zero polynomial. This means
that the images of these forms are a bases of $\fieldk[U,V]_d$ for
generic choice of $(a_j,b_j)$. Therefore we have on a generic line
$$\Syz(f_1, \ldots,f_{d+1})(d+1)|_{\PP^1}
\cong \Syz (U^d,U^{d-1}V \komdots V^d)(d+1) \cong \O_{\PP^1}^{d}\, .$$
So the restriction is semistable and hence the bundle itself on the projective space is semistable.

A generic set of $d+1$ forms of degree $d$ is generic on a generic line.
\end{proof}

\begin{corollary}
\label{tightgenerate}
Let $\fieldk$ denote an algebraically closed
field of characteristic $0$. Let $f_1, \ldots , f_{d+1} \in
\fieldk[X_0 \komdots X_N]$ be $d+1$ forms of degree $d$. Suppose
that there exists a linear morphism $\fieldk[X_0 \komdots
X_N] \ra \fieldk[U,V]$ such that the images of these forms are
linearly independent in $\fieldk[U,V]_d$. Then on a generic complete
intersection ring $R=\fieldk[X_0 \komdots X_N]/(G_1 \komdots
G_{N-1})$, where $G_j$ are generic forms of sufficiently high
degree, we have
$$(f_1 \comdots f_{d+1})^\star = (f_1 \comdots f_{d+1}) +R_{\geq d+1} \, .$$
This holds in particular for $d+1$ generic forms of degree $d$.
\end{corollary}
\begin{proof}
From Proposition \ref{restrictionline} and the Restriction Theorems
\ref{mehtaramanathanrestriction} or \ref{flennerrestriction}
it follows that the syzygy bundle is semistable on the smooth projective complete intersection curve defined
by $(G_1 \komdots G_{N-1})$ for generic $G_j$ of sufficiently high degree.
Hence the numerical formula from the introduction holds with the degree bound
$\sum_{i=1}^{d+1} \deg(f_i)/d = d(d+1)/d=d+1$.
\end{proof}

\begin{example}
The preceding Proposition and Corollary are applicable for $d+1$
forms of type
$$X^d+ ZP_0(X,Y,Z),\, X^{d-1}Y+ ZP_1(X,Y,Z),\, \ldots, \,Y^d+ ZP_{d}(X,Y,Z),\,$$
where $P_i$ are polynomials of degree $d-1$. By setting $Z=0$, these
forms yield all monomials of $\fieldk[X,Y]_d$.

The easiest instance is given by setting $P_i=0$, which gives just
all monomials in $X$ and $Y$. Here the equality $\Syz(X^d \komdots
Y^d)(d+1) \cong  \O^{d}$ holds already on $D_+(X,Y) \subset \PP^2$
and the stated result is true for every curve $V_+(G) \subset
D_+(X,Y)$ (this condition means that $X$ and $Y$ are parameters in
$R=\fieldk[X,Y,Z]/(G)$). An element $f$ of degree $m$ yields a
cohomology class in $H^1(V_+(G),\O(m-d-1)^{d})$ and one may argue on
the components. In this special case the formula in Corollary
\ref{tightgenerate} is also clear from \cite[Theorem
5.11]{hunekesmithkodaira}.
\end{example}

\begin{example}
\label{xyz21}
Consider $X^3,Y^3,Z^3,X^2Y$. Setting $Z=X+Y$, the restriction
yields four independent polynomials.
Hence the bundle is
semistable and it follows that
$R_{\geq 4} \subseteq (X^3,Y^3,Z^3, X^2Y)^\star$
in a generic hypersurface ring
$\fieldk[X,Y,Z]/(G)$.
\end{example}

\begin{question}
Let  $n$ monomials in $ \fieldk[X_0 \komdots X_N]$ of the same
degree $d$ ($n \leq d+1$) be given. When does there exist a linear
mapping
$$\fieldk[X_0 \komdots X_N] \lra \fieldk[U,V]$$
such that the images of the monomials are linearly independent?
\end{question}

\begin{example}
Consider the five monomials
$$X^4,Y^4,Z^4,X^3Y,X^3Z \in \fieldk[X,Y,Z]$$
of degree four and let $\shS = \Syz(X^4,Y^4,Z^4,X^3Y,X^3Z )$ denote
their syzygy bundle. The images of the monomials $X^4, X^3Y$ and
$X^3Z$ are linearly dependent for every linear homomorphism
$\fieldk[X,Y,Z] \ra \fieldk[U,V]$. It follows that for every line
$\PP^1 \subset \PP^2$ the restriction $\shS|_{\PP^1}$ is not
semistable, since the dependence yields non-trivial sections in
$(\shS|_{\PP^1})(4)$
(but the degree of $\shS(4)$ is $-4$).

The syzygy bundle on $\PP^2$ has global sections of degree $5$, and $\Gamma(\PP^2, \shS(5))$ is spanned by
$$(Z,0,0,0,-X) ,\,  (Y,0,0,-X,0) \mbox{ and } (0,0,0, Z,-Y)\, .$$
These global syzygies span a subsheaf which is isomorphic to $\Syz(X,Y,Z)(2)$.
Its slope is $1/2$, whereas $\shS(5)$ has slope $0$, hence $\shS$ is not semistable.

We want to compute its Harder-Narasimhan filtration.
We have the exact sequence
$$0 \lra \Syz(X,Y,Z)(m-3) \lra \shS(m)\stackrel{(p_2,p_3)}{ \lra} \O(m-4) \oplus \O(m-4) \, .$$
The syzygy subbundle on the left is semistable. The image of the
last mapping is a torsion-free subsheaf of rank $2$. This quotient
is given locally by
$$ \shQ(m)= \{ (s,t) \in \O(m-4) \oplus \O(m-4) : sY^4+tZ^4 \in (X^3) \}  \, .$$
For $m$ large enough, $\shQ(m)$ is generated by its global sections,
hence we look at this condition on $\fieldk[X,Y,Z]$. Then either
both $s$ and $t$ are multiples of $X^3$ or $sY^4+tZ^4=0$. This gives
the resolution
$$  0 \lra \O(m-11) \lra    \O(m-7) \oplus \O(m-7) \oplus \O(m-8)      \lra \shQ(m) \lra 0 \, ,$$
where the surjection is given by $1 \mapsto (X^3,0)$, $1 \mapsto (0, X^3)$ and $1 \mapsto (Z^4,-Y^4)$
and the injection is given by $1 \mapsto ( -Z^4,Y^4 ,X^3)$.
The quotient sheaf has degree $\deg Q(m)= 2m-11$ and it is semistable, since its first non-trivial section is for
$m=7$.
So we have found the Harder-Narasimhan filtration of $\shS=\Syz(X^4,Y^4,Z^4,X^3Y,X^3Z)$.
\end{example}

\section{Wedge criteria for stability}
\label{wedgesection}

Let $\shS $ denote a coherent torsion-free sheaf on $\PP^N$. A
coherent subsheaf $\shT \subseteq \shS$ of rank $r$ yields
$\bigwedge^r \shT \ra \bigwedge^r \shS$. The bidual of $\bigwedge^r
\shT$ is an invertible sheaf and its degree is by definition the
degree of $\shT$. Therefore the maximal degree of a subbundle of
rank $r$ is related to the global section $\neq 0$ of $(\bigwedge^r
\shS)(k)$. In particular we have the following criterion for
semistability, see \cite[Proposition 1.1 and the following remark
there]{bohnhorstspindler}.

\begin{proposition}
Let $\shS $ denote a locally free sheaf on $\PP^N$ over an
algebraically closed field $\fieldk$ of characteristic $0$. Then
$\shS$ is semistable if and only if for every $r < \rk (\shS)$ and
every $ k < -r \mu(\shS)$ there does not exist a global section
$\neq 0$ of $(\bigwedge^r \shS) (k)$.
\end{proposition}

\begin{proof}
If $\shS$ is semistable, then all its exterior powers
$\bigwedge^r \shS$ are also semistable (in characteristic $0$)
due to \cite[Corollary 3.2.10]{huybrechtslehn}.
Hence $(\bigwedge^r \shS ) \otimes \O(k) $ does not have global sections $\neq 0$
for $\mu ((\bigwedge^r \shS ) \otimes \O(k)) < 0$, which means that
$k < - \mu ( \bigwedge^r \shS ) =-r \mu( \shS)$.

Now suppose that the condition on the global sections is fulfilled
(this direction holds in any characteristic), and let $\shT \subset
\shS$ denote a coherent subsheaf of rank $r$. Then $\bigwedge^r \shT
\subset \bigwedge^r \shS$ and $(\bigwedge^r \shT) ^{\dual \dual}
\cong \O(m)$ is an invertible sheaf on $\PP^N$, where $m= \deg
(\shT)$. But then also $ \bigwedge^r \shT \cong \O(m)$ outside a
closed subset of codimension $\geq 2$. Since $\shS$ is locally free,
this gives a non-trivial morphism $\O(m) \ra \bigwedge^r \shS$.
Therefore $(\bigwedge^r\shS)(-m)$ has a global section $\neq 0$, so
$ -m \geq - r \mu (\shS)$ by assumption and hence $\mu(\shT)=
\frac{m}{r} \leq \mu(\shS)$. \end{proof}

So if we want to apply this stability criterion we need to get control on the exterior powers of
a syzygy bundle
$\Syz(f_i, i \in I)$ and their global sections. From the embedding
$$ \Syz(f_i,\, i \in I) \hookrightarrow \bigoplus_{i \in I} \O(-d_i)$$
we get the canonical embedding
$$ \bigwedge^r (\Syz(f_i,\, i \in I)) \lra  \bigwedge^r ( \bigoplus_{i \in I} \O(-d_i))
\cong  \bigoplus_{J \subseteq I,\, |J|=r } \O( - \sum _{i \in J} d_i) \, .$$
Here the identification on the right
is given by sending
$$s_1 \wedgedots s_r  \longmapsto \det \big( ( s_{ji} )_{j=1 \komdots r,\, i \in J} \big) \, .$$

\begin{lemma}
\label{wedgedescribe}
Let $\fieldk$ denote a field. Let $f_i \in
R=\fieldk[X_0 \komdots X_N]$, $i \in I=\{1 \comdots n\}$, denote
homogeneous, $R_+$-primary polynomials.
Then the sequence
$$0 \lra \bigwedge^r (\Syz(f_i, i \in I))
\lra \bigoplus_{|J|=r} \O(- \sum_{i \in J} d_i)
\stackrel{\varphi}{\lra} \bigoplus_{|K|=r-1} \O(- \sum_{i \in K} d_i) $$
is exact on $\PP^N$, where $\varphi$ sends
$ e_J \longmapsto  \sum_{k \in J} \sign(k,J) f_k e_{J- \{k\}}$
{\rm(}we use the induced order on $J \subseteq I$ to define
$\sign(k,J)$
\ifthenelse{\boolean{paper}}{to be $1$ if $k$ is an even element in $J$ and $-1$ otherwise}{
as in Lemma \ref{detidentify}}{\rm)}.
\end{lemma}
\begin{proof}
This is a global version of the local fact that $\bigwedge^r(V
\oplus R) \cong \bigwedge^r V \oplus \bigwedge^{r-1} V$ for a free
$R$-module $V$.
\ifthenelse{\boolean{paper}}{}{
We write down the sequence locally
on $D_+(f_1)$. We have the identification
$$(\O(-d_2) \oplusdots   \O(-d_n))|_{D_+(f_1)} \cong \Syz(f_i, i \in I) |D_+(f_1) $$
given by
$$(a_2 \komdots a_n) \longmapsto (- \frac{\sum_{i=2}^n a_if_i}{f_1} ,a_2 \komdots a_n) \, .$$
This identification yields the identification
$$\bigwedge ^r (\Syz(f_i, i \in I)) |D_+(f_1)
\cong \bigoplus_{|J|=r, 1 \not \in J} \O(- \sum_{i \in J}
d_i)|D_+(f_1) \, .$$ For a subset $J \subseteq I$, $|J|=r$, $1
\not\in J$ the rational $r$-form
$$\bigwedge_{i \in J}(- \frac{1}{f_1}, 0 \komdots 0, \frac{1}{f_i},0 \komdots 0)$$
corresponds under this identification to the section $\prod_{ i \in
J} \frac{1}{f_i}$ of $\O(- \sum_{i \in J} d_i)$ in the $J$-th
component and to $0$ in the other components. Therefore the first
mapping in the sequence is (under this identification) given by
$$e_J \longmapsto e_J + \sum_{k \in J} \sign(k,J)  \frac{f_k}{f_1} e_{\{1\} \cup J -\{k \}}\, .$$
The composition of the first mapping with $\varphi$ gives then
\begin{eqnarray*}
& &
\varphi(e_J + \sum_{k \in J} \sign (k,J)  \frac{f_k}{f_1} e_{\{1\} \cup J -\{k \}} ) \cr
&=&
\sum_{k \in J} \sign (k,J) f_k e_{J- \{k\}}
 +\sum_{k \in J} \sign(k,J) \frac{f_k}{f_1} \varphi(e_{\{1\} \cup J- \{k\}})
\cr
&=&
\sum_{k \in J} \sign (k,J) f_k e_{J- \{k\}}  -\sum_{k \in J} \sign(k,J) \frac{f_k}{f_1} f_1 e_{ J- \{k\}}
\cr
& &
+\sum_{k \in J} \sign(k,J) \frac{f_k}{f_1}
(\sum_{j \neq k} \sign (j,\{1\} \cup J- \{k\})  f_j ) e_{ \{1 \} \cup J -\{k ,j \}} \cr
&=& \sum_{j \neq k} c(k,j)
\frac{f_kf_j}{f_1} e_{ \{1 \} \cup J -\{k ,j \}} \, ,
\end{eqnarray*}
$c(k,j)= \sign (k,J) \sign(j, \{1\} \cup J- \{k \}) + \sign(j,J) \sign(k, \{1\} \cup J - \{j\})$.
But these coefficients are $=0$, since for $k < j$ we have
$\sign (k,J) = -\sign (k, \{1\} \cup J- \{ j\}) $ and
$\sign( j,J)= \sign (j,\{1\} \cup J- \{k \})$.

Now suppose that $ \varphi$ sends $\sum_J a_J e_J$ to $0$. In the
image of the first mapping we have the term $e_J$ for $J$, $1
\not\in J$, and all other expressions do contain $e_K$ with $1 \in
K$. Hence we may assume by adding elements of the image that $a_J=0$
for all $J$ with $ 1 \not\in J$. The image of $a_Ke_K$ ($1 \in K$)
under $\varphi$ contains the expression $f_1 a_K e_{K-\{1\}}$, but
this component is reached by no other element. Therefore $a_K=0$.}
\end{proof}

\begin{remark}
With the results of this section it is in principal possible to
decide whether a given syzygy bundle $\Syz(f_1 \komdots f_n)$ is
semistable or not. The exterior bundles $\bigwedge ^r \Syz$ are
given as kernels of some mappings between splitting bundles, hence
the minimal degree of a global section $\neq 0$ is computable with
Groebner basis techniques. An algorithm for this was developed and implemented by
A. Kaid in the computer algebra system CoCoA.
\end{remark}

\begin{remark}
\label{degreeremark}
Let $\shS \hookrightarrow \bigoplus_{\indi \in I} \O(-d_\indi)$ be a subsheaf of rank
$r$. We describe a method to compute $\deg (\shS)$.
Let $s_1 \comdots s_r \in \Gamma(\PP^N, \shS(m))$ be $r$ global sections which are linearly independent in at least one point
(hence on an open subset).
This gives a mapping
$\O^r \longrightarrow \shS(m) \hookrightarrow \bigoplus_{\indi \in I} \O(m-d_\indi)$.
Let these sections be given as $s_\indsec =(\fug_{\indsec \indi}) $, $\fug_{\indsec \indi} \in
\Gamma(\PP^N,\O(m-d_\indi))$, $\indsec=1 \comdots r$.
For every $J \subseteq I$
we look at the projection
$\bigoplus_{\indi \in I} \O(-d_\indi) \to \bigoplus_{\indi \in J} \O(-d_\indi) $
and the induced mapping
$\O^r \longrightarrow \shS(m) \longrightarrow \bigoplus_{\indi \in J}
\O(m-d_\indi)$,
and at
$$\O \cong \bigwedge^r \O^r \longrightarrow \bigwedge^r (\shS(m)) \longrightarrow \bigwedge^r(\bigoplus_{\indi \in J} \O(m-d_\indi))
\cong \O(rm- \sum_{\indi \in J} d_\indi ) \, ,$$
which is given by
$1 \mapsto \det ( (\fug_{\indsec \indi})_{1 \leq \indsec \leq r, \indi \in J}) =: \fuh_J$.
These $\fuh_J$ give together a map
$\O \cong \bigwedge^r \O^r \longrightarrow \bigwedge^r (\shS(m))
\longrightarrow
\bigoplus_{J \subseteq I, |J|=r} \O(rm- \sum_{\indi \in J} d_\indi )
$.
The degree of  $\bigwedge^r
(\shS(m))$ can be computed by counting the zeroes of codimension one of this section in
$\bigwedge^r( \shS(m))$, and this can be estimated by counting the zeroes of codimension one of this section in the direct
sum, which is the degree of the highest common factor of all
$\fuh_J$,
$|J|=r$. So we get the estimate $\deg (\shS(m)) \leq
\deg(\hcf(\fuh_J, |J|=r))$.

Suppose now that we have a syzygy bundle.
Global sections of
$\Syz(f_\indi, \indi \in I)(m)$ (where $m:= \sum_{\indi \in I} d_\indi$) are sometimes
(see the proof of Theorem \ref{maximalslope})
given as
$$ s= ( a_{\indi} \prod_{\indgen \neq \indi} f_\indgen )_{\indi \in I} $$
with $\sum_{\indi \in I} a_\indi =0$.
Suppose that $r$ such sections $s_1 \comdots s_r$
are given and are global sections of a subsheaf $\shS(m)$ of
$\Syz(f_i)(m)$ of rank
$r$ which are linearly independent in a point. Write
$s_\indsec = ( a_{\indsec \indi} \prod_{\indgen \neq \indi} f_\indgen )_{\indi \in I} $.
Hence for a subset $J \subseteq I$ with $r$ elements we get
$\fuh_J
=\det (  (a_{\indsec  \indi} \prod_{\indgen \neq \indi}  f_\indgen)_{1 \leq  \indsec  \leq r,  \indi \in J})$,
and this is
(as the expressions
$\prod_{\indgen \neq \indi} f_\indgen$ are constant in the
column corresponding to $\indi \in J$)
$$
(\prod_{\indgen \in I} f_\indgen)^{r-1} (\prod_{\indgen \in I-J} f_\indgen)
 \cdot \det (( a_{\indsec \indi})_{1 \leq  \indsec  \leq r,  \indi \in J})
\, .
$$
So here the highest common factor of the expressions
$\prod_{\indgen \in I-J} f_\indgen$ for $|J|=r$ and $\det (( a_{\indsec
\indi})_{1 \leq  \indsec  \leq r,  \indi \in J}) \neq 0$ is crucial.
We get then
\begin{eqnarray*}
\deg \shS \!\!
&=& \! \deg (\shS(m)) - r m \cr
&\leq&\! (r-1)m +\deg (\hcf  ( \prod_{\indgen \in I-J} f_\indgen, | J|=r, \det (( a_{\indsec
\indi})_{1 \leq  \indsec  \leq r,  \indi \in J}) \neq 0  )) -rm
\cr &=& \!
- \sum_{\indgen \in I} d_\indgen +\deg (\hcf  ( \prod_{\indgen \in
I-J} f_\indgen, | J|=r, \det (( a_{\indsec
\indi })_{1 \leq  \indsec  \leq r,  \indi \in J} \neq 0  ))\, .
\end{eqnarray*}
\end{remark}

\section{Stability of syzygies bundles of monomial ideals}
\label{monomialsection}

We consider now the case where $f_i \in R=\fieldk[X_0 \komdots
X_N]$, $i \in I$, are monomials and we will write
$f_i = X^{\tupmon_i} = X_0^{\tupmon_{i 0}} \cdots X_N^{\tupmon_{i N}}$, where $\tupmon_i \geq 0$ are integral lattice points
in $\NN^{N+1}$. Their degrees are $d_i=| \tupmon_i|= \sum_{j=0}^N
\tupmon_{i j}$.
We will apply the theory of toric bundles which has been developed
by Klyachko (see \cite{klyachkoselecta},
\cite{klyachkospectralproblems},
\cite{klyachkoequivariant}). We consider the projective space
$\PP^N$ as a toric variety with the torus $T= G_m^N =(\AA^\times)^N$ acting as
$$(t_1 \comdots t_N)(\coox_0 \comdots \coox_N)=( \coox_0,t_1 \coox_1 \comdots t_N \coox_N) \, .$$
A toric bundle is a vector bundle $V \to \PP^N$
with an action of $T$ on $V$ compatible with the basic torus action.

For every tuple $\tup=(\tup_0 \comdots \tup_N)$ we can make the line bundle
$\O(\sum_j \tup_j) \to \PP^N$ into a toric line bundle by the action
$$(t_1 \comdots t_N)(\coox_0 \comdots \coox_N \seco \varex)= ( \coox_0,t_1 \coox_1 \comdots t_N \coox_N \seco t_1^{\tup_1} \cdots t_N
^{\tup_N} \varex) \, ,$$ which we denote by $\O(\tup)$.
Recall that $\O(k)$ ($k \in \ZZ$) itself (disregarding any toric structure)
is given by dividing through the $\AA^\times$-action
$\eles (\coox_0 \comdots \coox_N \seco \varex)=(\eles \coox_0 \comdots \eles \coox_N \seco \eles^{k} \varex)
$.
For a given family of monomials $X^{\tupmon_i}$
also $\bigoplus _{i \in I} \O(-\tupmon_i)$ is a toric bundle and
the monomials define a toric morphism of this sum to $\O$, hence
$\Syz(X^{\tupmon_i}, i
\in I)$ is a toric subbundle.
Explicitly, $(t_1 \comdots t_N)$ sends
a point $(\coox_0 \comdots \coox_N \seco \sect_i, i \in I)$
(with $\sum_{i \in I} \sect_i \coox^{\tupmon_i}=0$) over
$(\coox_0 \comdots \coox_N)$
to $(\coox_0 , t_1 \coox_1\comdots t_N \coox_N \seco
t_1^{-\tupmon_{i 1}}
\cdots t_N^{-\tupmon_{i N}} \sect_i  , i \in I)$.
Note that
$$\sum_{i \in I}   t_1^{-\tupmon_{i 1}} \cdots t_N^{-\tupmon_{i N}}
\sect_i \cdot
 \coox_0^{\tupmon_{i 0}} (t_1 \coox_1)^{\tupmon_{i 1}} \cdots  (t_N \coox_N)^{\tupmon_{i N}}
  =\sum_{i \in I} \sect_i x^{\tupmon_i}=0 \, .$$

Klyachko studies toric bundles $\buntor$ with the help of families of
filtrations in the special fiber $\fibspec_\point= \buntor \tensor
\kappa(\point)$, where
$\point$ is a closed point outside of any toric hypersurface. Every toric
hypersurface
$H_\alpha$ determines a decreasing filtration $\fibspec^\alpha (m)$,
$m \in
\ZZ$, of vector subspaces in $\fibspec_\point$.
For example, the toric line bundle $\O(\tup)$
on $\PP^N$ corresponds to the family of filtrations on $\fieldk$
given by $ \fieldk^\alpha (m)= \fieldk$ for $m \leq \tup_\alpha$ and
$ \fieldk^\alpha (m)= 0$ for $m > \tup_\alpha$ (compare \cite[Example 2.3]{klyachkoequivariant}).
The category of toric vector bundles
on a toric variety is equivalent to the category of vector spaces
with such families of filtrations fulfilling certain compatibility
conditions (see \cite[Theorem 2.2.1]{klyachkoequivariant}).

We collect some of the main properties which we need in the sequel
of this section.

\begin{lemma}
\label{toriclemma}
Let $\buntor$ be a toric bundle on $\PP^N_\fieldk$, $\fieldk$ an algebraically closed field, and set
$\fibspec = \fibspec_\point$ with the corresponding filtrations
$\fibspec^\alpha(m)$, $\alpha=0 \comdots N$,
where $\point$ is a closed point in the torus.
Then the following hold.
\numiii
\begin{enumerate}

\item
Let a vector $\vecspec \in \fibspec$ be given. Let $n_\alpha$ be the maximal integer
{\rm(}take $\infty$ for $\vecspec=0${\rm)}
such that
$\vecspec \in \fibspec^\alpha(n_\alpha)$. Then the $\fieldk$-linear mapping $\fieldk \to \fibspec$, $1 \mapsto \vecspec$,
extends to a toric bundle morphism
$\O(\tup) \to \buntor$ under the condition that $\tup_\alpha \leq n_\alpha$ for $\alpha= 0 \comdots N$.

\item
Let a linear form  $\foli: \fibspec \to \fieldk$ be given and let
$m_\alpha $ be the smallest number such that $\fibspec^\alpha(m_\alpha) \subseteq \ker \foli$
{\rm(}take $-\infty$ for $\foli=0${\rm )}.
Then $\foli$ extends to a toric bundle morphism $\buntor \to \O(\tup)$ for $\tup=(\tup_\alpha)$,
$\tup_\alpha \geq m_\alpha$.

\item
Let $ \subfibspec \subseteq \fibspec$ be a vector subspace.
If $\foli_\indk: \fibspec \to \fieldk$
is a family of linear forms such that $\subfibspec = \bigcap_\indk
\ker \foli_\indk$,
then the kernel sheaf {\rm(}which is not locally free in
general{\rm)} of the toric mapping
$( \folitor_\indk) :\buntor \to
\bigoplus_\indk
\O(m_\indk)$ {\rm(}as constructed in {\rm(ii)}{\rm)}
is a toric subsheaf $\subbuntor$ {\rm(}in the sense that it is a
toric subbundle over an open toric subvariety which contains all
points of codimension one{\rm )} such that
$\subbuntor
\tensor
\kappa(\point)= \subfibspec$. In particular,
subspaces of the special fiber correspond to toric subsheaves which
are locally free in codimension one.

\item
For $\vecspec \in \subfibspec \subseteq \fibspec$ the global morphism {\rm(}constructed in {\rm(i)}{\rm)} factors
through $\subbuntor$.

\item
The maximal destabilizing subsheaf of $\buntor$ is given by a toric subbundle
which is defined on an open toric subvariety containing all points
of codimension one.
\end{enumerate}
\end{lemma}
\begin{proof}
(i). Let $\varphi :\fieldk \to \fibspec$ be the map given by $1 \mapsto \vecspec$.
Then $\varphi(\fieldk^\alpha(m)) \subseteq \fibspec^\alpha(m)$ for all $m$
if and only if $\varphi(\fieldk^\alpha(\tup_\alpha) ) \subseteq
\fibspec^\alpha(\tup_\alpha)$
if and only if $\vecspec \in \fibspec^\alpha(\tup_\alpha)$ if and
only if $ \tup_\alpha \leq n_\alpha$ (this holds for every
$\alpha$).
Such a filtered linear mapping corresponds to a toric bundle
morphism
$\O(\tup) \to
\buntor$ by
Klyachko's theorem \cite[Theorem 2.2.1]{klyachkoequivariant}.

(ii).
The linear mapping $\foli: \fibspec \to \fieldk$ is for $\tup=(\tup_\alpha)$,
$\tup_\alpha \geq m_\alpha$, compatible with the filtrations, since for
$m < m_\alpha$ we have $m < \tup_\alpha$ and so $\fieldk^\alpha(m)=\fieldk$, and
for $m \geq m_\alpha$ we have $ \foli (\fibspec^\alpha(m))= 0$.
Hence again by \cite[Theorem 2.2.1]{klyachkoequivariant}
$L$ extends to a toric morphism $\buntor \to \O(\tup)$ for $\tup=(\tup_\alpha)$,
$\tup_\alpha \geq m_\alpha$.

(iii).
The linear forms $\foli_\indk: \fibspec \to \fieldk$ induce by (ii) a
toric morphism
$\buntor \to \bigoplus_\indk \O(m_\indk)$, where the $m_\indk= m_{\indk \alpha}$ are defined as in (ii).
The kernel sheaf is locally free in codimension one and its special
fiber is
$\subfibspec$.
Since the toric automorphisms respect the kernel it is a toric
subsheaf.
In particular, every subspace $\subfibspec \subset \fibspec$
is the special fiber of a toric subsheaf. On the other hand, let two
toric subsheaves $\subbuntor_1$ and $\subbuntor_2$ of $\buntor$ with
the same special fiber be given. Then the induced filtrations are
the same and so they are isomorphic as toric bundles on a certain
open toric subvariety. By \cite[Theorem 2.2.1]{klyachkoequivariant}
they must be the same subbundle, since the embedding is determined
by the filtered linear mapping.

(iv).
The composed mapping $\O(\tup) \to \buntor  \to \bigoplus_\indk
\O(m_\indk)$ is the zero map on the special fiber, since $\vecspec \in \subfibspec = \bigcap_\indk \ker \foli_\indk $.
Hence, again by \cite[Theorem 2.2.1]{klyachkoequivariant},
it is the zero map and so it factors through the kernel, which is
$\subbuntor$ by part (iii).

(v).
For $\PP^2$ this is proven in \cite[Theorem 3.2.2]{klyachkoselecta}, but in general we
have to be a bit more careful.
Let $\subbuntor$ be the maximal destabilizing subsheaf of $\buntor$. As this
is uniquely determined, we must have $t^*(\subbuntor) =
\subbuntor$ as a subsheaf of $\buntor$ for (the automorphism given
by) $t \in T$.
In particular,
$\shF$ is locally free on an open toric subvariety which contains all
points of codimension one. The action $t: \buntor \to \buntor$ must
send
$\subbuntor_\point$ to $\subbuntor_{t(\point)}$.
Hence the action restricts to $\subbuntor$ and so $\subbuntor$ is toric
(but not necessarily a bundle on the whole).
\end{proof}

\begin{lemma}
\label{linearlemma}
Let $v_\indi$, $\indi \in I$, be a set of spanning vectors $\neq 0$ in a vector space $U$ of
dimension $r $ and with $\sum_{\indi \in I} v_\indi=0$. Then there
exists a partition
$I=I_1 \uplusdots I_\numpart \uplus \tilde{I} $ such that
$ \sum_{\indpart=1}^\ell ( |I_\indpart| -1) =r$
and such that, setting
$V_\indpart= \langle v_\indi, \indi \in I_\indpart \rangle$, the following holds:
the set of vectors
$\{v_\indi : \indi \in I_\indpart\}$
is linearly dependent modulo the subspace $V_1 \plusdots
V_{\indpart-1}$,
but all strict subsets are independent.
\end{lemma}
\begin{proof}
Note first that for every hyperplane $H \subset U$ there exist at
least two vectors outside $H$, because of the sum property. We do
induction on
$r$.
For
$r=1$ any $I_1=\{i,j\}$ and $\tilde{I}= I-I_1$ will do.
So let $r \geq 2$. Reorder so that $v_1 \comdots v_r$ are
a basis of $V$.
Take $v_{r+1}$
(which must exist because of the sum property)
and consider
$\{v_1 \comdots v_r,v_{r+1}\}$.
If this set has the property that
every strict subset is independent, then we take $I_1=\{1
\comdots r+1\}$ ($\ell=1$) and we are done.
In the other case there exists a dependent
(strict) subset, which must contain $v_{r+1}$, since the first $r$ vectors
$v_\indi$ are independent. Then either this set has the property or we decrease the set further until we arrive at
a set with the required properties (the smallest possibility is that
of $\{v_\indi,v_{r+1}\}$ being dependent).

Let $I_0 \subseteq \{1 \comdots r,r+1\}$ be such an index set
and let $V_0 = \langle v_\indi, \indi \in I_0 \rangle$ ($\neq 0$) be
the subspace.
Let $I' \subset I$ be the subset consisting of the
$i$ such that $v_\indi$ do not belong to $V_0$.
Then the quotient space $V/V_0$ and the set of residue classes
$\{\bar{v}_\indi : \indi \in I' \}$ fulfill also all the assumptions
and is of smaller dimension.
Hence we apply the induction hypothesis
to get a partition
$I'=I'_1 \uplusdots I'_{\numpart} \uplus \tilde{I'}$ with the desired properties. Then
the sets
$I_0$,
$I_\indpart:= I'_\indpart $ ($\indpart =1 \comdots \ell$)
and
$\tilde{I} = \tilde{I'} \cup (I- (I_0 \cup I')) $
form a partition of $I$ with the desired properties.
\end{proof}

\begin{theorem}
\label{maximalslope}
Let $f_i
=X^{\tupmon_i}$, $i \in I $, denote a set of primary monomials in
$\fieldk[X_0 \komdots X_N]$ of degree $d_i=| \tupmon_i|$.
For $J \subseteq I$ denote by $d_J$ the degree of the highest common
factor of $f_i, i \in J$. Then the maximal slope of
$\Syz(f_i, i \in I)$ is
$$\mu_{\rm max} (\Syz(f_i, i \in I))
=  \max_{J \subseteq I, |J| \geq 2}  \{ \frac{d_J- \sum_{i \in J} d_i}{ |J|-1} \}  \, . $$
\end{theorem}
\begin{proof}
It is clear that $\geq$ holds.
By Lemma \ref{toriclemma}(v) we only have to consider toric subsheaves
$\subbuntor \subseteq \Syz(f_i, i \in I)$ (i.e., toric subbundles defined on an open toric subvariety containing
all points of codimension one).
These are in
one-to-one correspondence (Lemma \ref{toriclemma}(iii)) with
subspaces
$\subfibspec
\subseteq \fibspec$ inside the special fiber $\fibspec$ of the syzygy bundle over the point $\point=(1 \comdots 1)$
($\fibspec$ itself is the hypersurface in $\fieldk^n$ given by $\sum_{\indi=1}^{n} a_\indi =0$).
So let $\subfibspec \subseteq \fibspec$ be a subspace of dimension $r$, given by $r$
linearly independent vectors $ w_1 \comdots w_r$, where
$w_\indsec = \sum_{\indi=1}^n a_{ \indsec \indi} e_\indi$,
$\sum_{\indi=1}^{n} a_{\indsec \indi} =0$.
We look at the global sections
$$s_\indsec = ( a_{ \indsec \indi} \prod_{\indk \in I, \, \indk \neq \indi}
f_\indk)_{ \indi \in I}
 \in \Gamma(\PP^N,\Syz(f_\indi, \indi \in I)(m))
\, ,$$
(where $m=\sum_{\indi=1}^n d_\indi$), which have $w_\indsec$ as their values at $(1 \comdots 1)$.
These sections are toric sections (where $\Syz(m)$ has the natural toric structure induced by
$\oplus_{i \in I} \O( \sum_{\indk \neq \indi} \tupmon_\indk)$),
hence they coincide (up to the twist) with the section constructed
in Lemma \ref{toriclemma}(i).
These sections factor through the toric subsheaf $\subbuntor (m)$ (Lemma
\ref{toriclemma}(iv)).

By Remark
\ref{degreeremark}
we get the estimate
$$\deg(\subbuntor) \leq -\sum_{\indgen \in I} d_\indgen +
\deg (\hcf (  \prod_{\indgen \in I-J} f_\indgen, | J|=r, \det
(( a_{\indsec \indi})_{1 \leq \indsec \leq r,\indi \in J}) \neq 0)) \, . $$
Set $v_{\indi} =(a_{\indsec \indi})$, $\indi \in I$, considered in
the vector space $\fieldk^r$. Note that
$\sum_{\indi \in I} v_\indi =0$. By Lemma \ref{linearlemma} there exists a
partition
$I=I_1 \uplusdots I_\numpart \uplus \tilde{I} $ such that
$ \sum_{\indpart=1}^\ell ( |I_\indpart| -1) =r$
and such that, setting
$V_\indpart= \langle v_i, i \in I_\indpart \rangle$, the following holds:
the set
$\{v_i : i \in I_\indpart\}$
is linearly dependent modulo the subspace $V_1 \plusdots
V_{\indpart-1}$,
but all strict subsets are independent.
Then for all subsets
$$J=J_1 \uplusdots J_\numpart, \, J_\indpart \subset I_\indpart, |J_\indpart| = |I_\indpart| -1$$
the vectors $v_\indi$, $\indi \in J$, are linearly independent and
so the determinantial coefficients (for these $J$) are
$\neq 0$.
Hence
$$
\deg (\hcf  ( \!\! \prod_{\indgen \in I-J} \!\! f_\indgen, |J| \! =  r, \det
(\!( a_{\indsec \indi})_{1 \leq \indsec \leq r,\indi \in J}) \neq 0\!) \!)
\! \leq \!
\deg (\hcf (\!\!  \prod_{\indgen \in I-J}  \!\! f_\indgen, J \mbox{ as above})\!) .
$$
The products on the right can be written as
$(\prod_{\indi \in \tilde{I}} f_\indi)
f_{\indj_1} \cdots f_{\indj_\ell}$ for any choice $ \indj_1 \in
I_1
\comdots
\indj_\ell \in I_\ell $.
So their highest common factor is
$(\prod_{\indi \in \tilde{I}} f_\indi) \cdot \hcf (f_i, i \in
I_1) \cdots \hcf (f_i, i \in I_\ell)$.
Therefore
\begin{eqnarray*}
\deg(\subbuntor) &\leq &
-  \sum_{\indi \in I_1 \cupdots I_\ell}  d_i + \deg(\hcf (f_i, i \in
I_1))  \plusdots  \deg( \hcf (f_i, i
\in I_\ell)) \cr
&=& \sum_{\indpart=1}^\ell (\sum_{\indi \in I_\indpart}- d_\indi
+ \deg(\hcf(f_\indi, \indi \in I_\indpart))) \, .
\end{eqnarray*}
On the right we have the degree of the subsheaf
$\Syz(f_\indi, \indi \in I_1) \oplusdots \Syz(f_\indi , \indi \in I_\ell)$
of rank $r$. Its slope can be estimated by the maximum of the slopes
of its direct summands, which are
$\frac{d_{I_\indpart} -\sum_{i \in
I_\indpart} d_i}{ |I_\indpart|-1}$.
\end{proof}

We can now state our combinatorial criterion for a monomial family
to have a semistable syzygy bundle (the necessity of the condition
was already established in Proposition \ref{subsetcondition}).

\begin{corollary}
\label{monomstablecrit} Let $f_i =X^{\tupmon_i}$, $i \in I,$ denote a
set of primary monomials in $\fieldk[X_0 \komdots X_N]$ of degree
$d_i= |\tupmon_i|$. Suppose that for every subset $J \subseteq I$,
$|J| \geq 2$, the inequality
$$ \frac{ d_J - \sum _{i \in J} d_i}{|J|-1} \leq \frac{-\sum_{i \in I} d_i}{|I|-1}$$
holds, where $d_J$ is the degree of the highest common factor of
$f_i, i \in J$. Then the syzygy bundle $\Syz(f_i, i \in I)$ is
semistable {\rm(}and stable if $<$ holds{\rm)}.
\end{corollary}
\begin{proof}
This follows at once from
Theorem \ref{maximalslope}.
\end{proof}

\begin{corollary}
\label{monomialstablecrit}
Let $f_i =X^{\tupmon_i}$, $i =1 \komdots n
$ denote a set of primary monomials of the same degree $d$ in
$\fieldk[X_0 \komdots X_N]$. For every monomial $X^\tupsub$ of degree
$e=|\tupsub| \leq d$ let $s_{\tupsub}$ denote the number of monomials in
the family which are multiples of $X^{\tupsub}$. Then the syzygy
bundle
$\Syz(f_1 \komdots f_n)$ is semistable if and only if for every
$\tupsub$ the following inequality holds:
$$\frac{s_{\tupsub}-1}{d-e} \leq  \frac{n-1}{d} \, . $$
\end{corollary}
\begin{proof}
Let $J \subseteq I=\{1 \komdots n \}$
denote the subset of monomials which are multiples of $X^\tupsub$.
The numerical semistability condition is that (setting $e=|\tupsub|$, $s=|J|$)
$$\mu(\Syz(f_i, i \in J))= \frac{e-sd}{s-1} \leq \frac{-nd}{n-1}
= \mu(\Syz(f_i,i\in I))\, .$$
This is equivalent with $e(n-1)-sd(n-1) \leq -(s-1)nd$
and hence with $sd \leq nd-e(n-1)$ and with $(s-1)d \leq
(n-1)d-e(n-1)= (n-1)(d-e)$, so the result follows.
\end{proof}

\section{Examples of monomial ideals}
\label{examples}

We first deduce the following result of Flenner (see
\cite[Corollary]{flennerrestriction} and \cite[Corollary
6.5]{ballicorestriction}) from our numerical criterion.

\begin{corollary}
\label{fullmonomial}
Let $\fieldk$ denote a field. Then the syzygy bundle
of the family of all monomials $\in \fieldk[X_0 \komdots X_N]$ of
fixed degree $d$ is semistable.
\end{corollary}
\begin{proof}
We want to apply Corollary \ref{monomialstablecrit}, so let $X^\nu$
be a monomial of degree $|\nu|=e \leq d$. Every monomial of degree
$d-e$ gives a multiple of $X^\nu$ of degree $d$, so we have to show
that
$$\frac{\binom{N+d-e}{N} -1}{d-e} \leq \frac{\binom{N+d}{N}-1}{d} \, .$$
We may assume successively that $e=1$, so we have to show that
$ (\binom{N+d}{N} -1)(d-1) \geq d (\binom{N+d-1}{N}-1)$, which is
equivalent with
$ \binom{N+d}{N} (d-1) \! \geq \! d \binom{N+d-1}{N}-1$. This is true for
$N=1$ or $d=1$ and follows for $N,d \geq 2$ from
$(N+d)(d-1) \geq d^2$.
\end{proof}

For a family consisting only of some powers of the variables we have
the following result, which is also a special case of Corollary
\ref{bohnhorstcorollary} and follows also from the Theorem
\ref{bohnhorstspindler} of Bohnhorst-Spindler.

\begin{corollary}
\label{powers} Consider the family $X_i^{d_i}$, $i=0 \komdots N$ in
$\fieldk[X_0 \komdots X_N]$. Suppose that $1 \leq d_0 \leqdots d_N$.
Then the syzygy bundle $\Syz(X_0^{d_0} \komdots X_N^{d_N})$ is
semistable on $\PP^N$ if and only if $(N-1) d_N \leq
\sum_{j=0}^{N-1} d_j $ holds.
\end{corollary}
\begin{proof} The numerical condition is necessary due to Corollary
\ref{numcondition}. On the other hand, again due to Corollary
\ref{numcondition} the necessary numerical conditions for smaller
ranks are also fulfilled, so the result follows from Corollary
\ref{monomstablecrit}. \end{proof}

We give some examples of small families of monomials in three
variables and check whether their syzygy bundles are stable or not.

\begin{corollary}
\label{fourmonomials} Let $X^{d_1}, Y^{d_2}, Z^{d_3}$ and $
X^{a_1}Y^{a_2}Z^{a_3}$ be four monomials in $\fieldk[X,Y,Z]$, $a_j <
d_j$. Set $d_4=a_1+a_2+a_3$. Then the syzygy bundle is semistable if
and only if the following two numerical conditions hold:

\renewcommand{\labelenumi}{(\roman{enumi})}

\begin{enumerate}

\item
$3 \max(d_1,d_2,d_3,d_4) \leq d_1 +d_2 + d_3+d_4$

\item
$\min (a_1+a_2+d_3, a_1+d_2+a_3, d_1+a_2+a_3, d_1+d_2, d_1+d_3, d_2+d_3) \geq \frac{\sum_{i=1}^4 d_i}{3}$.

\end{enumerate}
\end{corollary}
\begin{proof}
We apply the semistability criterion Corollary \ref{monomstablecrit}
for subsets $J$ with $|J|=2$ or $3$.
For $|J|=3$ we have $d_J=0$, so the condition is that
$\frac{- \sum_{i \in J} d_i}{2}  \leq \frac{ - \sum _{i=1}^4 d_i}{3}$, which is equivalent with
$- \sum_{ i \neq k} d_i \leq -2d_k$
for every $k$, so this is condition (i).
For $|J|=2$
the condition is
$d_J - \sum_{i \in J} d_i \leq \frac{- \sum _{i=1}^4
d_i}{3}$ or
$ \sum _{i \in J} d_i -d_J \geq \frac{\sum_{i=1}^4
d_i}{3}$ for all subsets $J$, $|J|=2$, which is condition (ii).
\end{proof}

\begin{example}
Consider $X^4,Y^4,Z^4,XYZ^2$. The first condition is clearly
satisfied. The minimum in the second condition is $6$ which is $\geq
16/3$, so the syzygy bundle is semistable.
If we replace however $XYZ^2$ by $XZ^3$, then the first condition is
again satisfied, but the minimum in the second condition is now $5$
and so this syzygy bundle is not semistable.
For $X^3,Y^3,Z^3,XY^2Z^2$ the second
condition is fulfilled, but the first is not fulfilled.
\end{example}

We consider now examples of more than four monomials.

\begin{example}
Consider now the monomials $X^6, Y^6, Z^6, X^2Y^2Z^2, XY^2Z^3$.
Their syzygy bundle is not semistable, since its slope is
$-\frac{30}{4} =-7.5$, but the subbundle $\Syz(X^2Y^2Z^2, XY^2Z^3)$
has slope $5-12 =-7$. This is also the maximal slope of this bundle.

For the monomials $X^6, Y^6, Z^6, X^2Y^2Z^2, X^3Z^3$ the slope
is again $-\frac{30}{4} = -7.5$. For $|J|=2$ the highest common
factor has maximal degree $4$, which gives slope $4-12=-8$. For
$|J|=3$ the common factor has degree at most $2$, which gives the
slope $\frac{2-18}{2} =-8$, so this syzygy bundle is stable.
\end{example}

\begin{example}
Consider the monomial family given by $X^5,X^4Z,Y^5,Y^4Z,Z^5 \!$, so
that the slope is $-6.25$. The subsheaf
$$\Syz(X^5,X^4Z) \oplus \Syz(Y^5,Y^4Z) \subset \Syz(X^5,X^4Z,Y^5,Y^4Z,Z^5)$$
has slope $-6 > -6.25 $ and it is the maximal destabilizing
subsheaf.
The subsheaves given by a subfamily of three elements do not
contradict semistability.
\end{example}

\begin{example}
\label{minimalsectionnotsubsheaf}
For a fixed $r$ the minimal degree of a global section of
$\bigwedge^r \Syz(f_i, i \in I)$ does in general not arise from a
subsheaf of rank $r$ of the form $\bigoplus_{\indpart =1}^\ell
\Syz(f_i, i \in I_\indpart)$
(though the maximal slope can be computed using only these subsheaves).
Look at the example given by the six monomials
$$X^4Y^2,\, X^4Z^2,\, Y^3Z^3,\, Y^5,\, Z^5,\,X^7 \, .$$
Their syzygy bundle is semistable according to the monomial
criterion (but not stable). For $r=2$, the subfamilies of three
elements yield global sections of
$(\bigwedge^2 \Syz(f_i, i \in I))(15)$,
but not of smaller degree, and the subsheaves of form
$\Syz(f_i,f_j) \oplus \Syz(f_s,f_t) $
yield only global sections of degree $\geq 16$.

There exists however also a section of degree $14$ of
$\bigwedge^2 \Syz(f_i, i \in I)$.
The subfamily $(X^4Y^2,X^4Z^2,Y^3Z^3)$ yields the section of degree
$18$
$$s_1 = -Y^3Z^3 e_{\{1,2\}} + X^4Z^2e_{\{1,3\}}  - X^4Y^2 e_{\{2,3\}} $$
(given in terms of Lemma \ref{wedgedescribe})
and $\Syz(X^4Y^2,Y^5) \oplus \Syz(X^4Z^2,Z^5)$ yields the section
$$s_2= Y^3Z^3 e_{\{1,2\}}
-X^4Y^3 e_{\{1,5\}} + X^4Z^3 e_{\{2,4\}} + X^8 e_{\{4,5\}}
\, .$$
Then $s_1+s_2$ is a multiple of $X^4$ and yields the section of degree $14$
$$ Z^2e_{\{1,3\}} -Y^3 e_{\{1,5\}} - Y^2 e_{\{2,3\}}
  +Z^3 e_{\{2,4\}} +
X^4 e_{\{4,5\}}
\, .$$
\end{example}

\begin{question}
Does there exist for every $d$ and every $n \leq \binom{N+d}{N}$ a
family of $n$ monomials in $\fieldk[X_0 \komdots X_N]$ of degree $d$
such that their syzygy bundle is semistable? This is due to
Corollary \ref{monomstablecrit} a purely combinatorial problem. A
positive answer to this question would imply that also the syzygy
bundle for generic polynomials of constant degree is semistable.
For $N=1$ (two variables) this is clearly true for the family
$X_0^d, X_0^{d-1}X_1 \komdots  X_0^{d-n+1}X_1^{n-1}$. In three
variables this is proved in \cite{brennermonomial}.
\end{question}

\section{Syzygy bundles of generic forms}
\label{generalsection}

What can we say about stability properties of $\Syz(f_1 \komdots
f_n)$ for generic homogeneous forms $f_1 \komdots f_n \in
\fieldk[X_0 \komdots X_N]$ of given degrees $d_i$? There is no hope
for semistable syzygy bundles unless the degrees satisfy the
necessary numerical condition described in Proposition
\ref{numcondition}. On the other hand, if these numerical conditions
are fulfilled, e. g. if the degrees are constant, then it is not
clear at all whether there exist semistable syzygy bundle of this
degree type. The degrees determine the Chern classes of the syzygy
bundle and therefore the question is equivalent to the following.
Does the moduli space $\shM(n-1, c_j)$ of rank $n-1$ stable vector
bundles on $\PP^N$ with Chern classes $c_j$ contain syzygy bundles?

We will give here some partial results for semistability using
results of Bohnhorst and Spindler \cite{bohnhorstspindler} and of
Hein (see the appendix).

\begin{theorem}{\rm(}Bohnhorst-Spindler{\rm )}
\label{bohnhorstspindler}
Let $\shE$ denote a vector bundle of rank $N$ on
the projective space $\PP^N$ over an algebraically closed field of characteristic $0$.
Suppose that $\shE$ has a resolution
$$0 \lra \bigoplus _{i=1}^k \O(a_i) \lra \bigoplus _{j=1}^{N+k} \O(b_j) \lra \shE \lra 0 \, . $$
and suppose that the pair $(a,b)$ is admissible, that means that
$a_1 \geqdots a_k$, $b_1 \geqdots b_{N+k}$ and
$a_i < b_{N+i}$ for $i=1 \komdots k$.
Then the following are equivalent.

\renewcommand{\labelenumi}{(\roman{enumi})}
\begin{enumerate}

\item
$\shE$ is semistable.

\item
$b_1 \leq \mu(\shE)= \frac{1}{N} (\sum_{j=1}^{N+k} b_j - \sum_{i=1}^k a_i
)$.

\end{enumerate}

\end{theorem}
\begin{proof}
See \cite[Theorem 2.7]{bohnhorstspindler}.
\end{proof}

\begin{corollary}
\label{bohnhorstcorollary}
Let $\fieldk$ denote an algebraically
closed field of characteristic $0$. Suppose that $f_1 \komdots
f_{N+1} \in \fieldk[X_0 \komdots X_N]$ are homogeneous polynomials
of degree $ d_1 \geqdots d_{N+1} \geq 1$. Suppose that $d_1 \leq
\frac{\sum_{i=2}^{N+1} d_i}{N-1}$ and that the $f_i$ are parameters.
Then the syzygy bundle $\Syz(f_1 \komdots f_{N+1})$ is semistable.
\end{corollary}
\begin{proof}
Since the $f_i$ are parameters their syzygy bundle is locally free and the presenting sequence
$$0 \lra \Syz(f_1 \komdots f_{N+1}) \lra \bigoplus_{j=1}^{N+1} \O(-d_j) \lra \O \lra 0$$
is exact on the right. Its dual is then also exact and we are in the
situation of the Theorem of Bohnhorst-Spindler with $k=1$, $a_1 =0$
and $b_j =d_j$. This pair is clearly admissible. The numerical
condition in the assumption is equivalent to the numerical condition
in Theorem \ref{bohnhorstspindler}(ii). Hence $\Syz^\dual$ is
semistable and then by \cite[Lemma
II.1.2.4]{okonekschneiderspindler} also $\Syz (f_1 \komdots
f_{N+1})$ is semistable. \end{proof}

\begin{corollary}
\label{parameter} Let $\fieldk$ denote an algebraically closed field
of characteristic $0$. Suppose that $f_1 \komdots f_{N+1} \in
\fieldk[X_0 \komdots X_N]$ are homogeneous parameters with degrees
$d_1 \geqdots d_{N+1} \geq 1$ such that $d_1 \leq
\frac{\sum_{i=2}^{N+1} d_i}{N-1}$. Then for generic forms  $G_1
\komdots G_{N-1} \in \fieldk[X_0 \komdots X_N]$ of sufficiently high
degree the equation
$$(f_1 \komdots f_{N+1})^\star = (f_1 \komdots f_{N+1}) + R_{\geq \frac{\sum_{i=1}^{N+1} d_i}{N} }$$
holds in $R=\fieldk[X_0 \komdots X_N]/(G_1 \komdots G_{N-1})$.
\end{corollary}
\begin{proof}
This follows from Corollary \ref{bohnhorstcorollary}, the restriction theorems and the
numerical formula for tight closure.
\end{proof}

\begin{corollary}
\label{genericsamedegree}
Let $\fieldk$ denote an algebraically
closed field of characteristic $0$. Let $f_1 \komdots f_{N+1} \in
\fieldk[X_0 \komdots X_N]$ denote $N+1$ generic homogeneous
polynomials of the same degree $d \geq 1$. Let $G_1 \komdots G_{N-1}
\in \fieldk[X_0 \komdots X_N]$ denote generic forms of sufficiently
high degree. Then
$$(f_1 \komdots f_{N+1})^\star = (f_1 \komdots f_{N+1}) + R_{\geq \frac{(N+1)d}{N} }$$
holds in $R=\fieldk[X_0 \komdots X_N]/(G_1 \komdots G_{N-1})$.
\end{corollary}

\begin{proof} $N+1$ generic homogeneous elements are parameters in
$\fieldk[X_0 \komdots X_N]$, so this follows from Corollary
\ref{parameter}. \end{proof}

\begin{remark}
Corollaries \ref{bohnhorstcorollary} and \ref{parameter}
generalize the case $N=2$ treated in Corollaries \ref{threeregularsemistable}
and \ref{threeprimarydegree}.
Corollary \ref{powers} is also a special case of Corollary \ref{bohnhorstcorollary}.
On the other hand,
we may deduce Corollary \ref{genericsamedegree} from Corollary \ref{powers}
without using the result of Bohnhorst-Spindler:
since semistability is an open property in a flat family
it is enough to establish the semistability property for a single choice
of homogeneous forms with given degree.
\end{remark}

\begin{theorem}{\rm(}Hein{\rm)}
\label{hein} Let $\fieldk$ denote an algebraically closed field and
let $f_1 \komdots f_n \in \fieldk[X_0 \komdots X_N]$, $N \geq 2$,
denote generic homogeneous polynomials of the same degree $d$.
Suppose that $n \leq d (N +1)$. Then their syzygy bundle $\Syz(f_1
\komdots f_n)$ is semistable on $\PP^N$.
\end{theorem}
\begin{proof}
See Theorem A.1 of the Appendix by G. Hein.
\end{proof}

\begin{corollary}
\label{genericgeneric} Let $\fieldk$ denote an algebraically closed
field of characteristic $0$. Let $f_1 \komdots f_n \in \fieldk[X_0
\komdots X_N]$ denote generic homogeneous forms of degree $d$, $n
\leq d (N+1)$. Then for a generic complete intersection ring
$R=\fieldk[X_0 \komdots X_N]/(G_1 \komdots G_{N-1})$ of sufficiently
high degree we have
$$ (f_1 \comdots f_n)^\star = (f_1 \comdots f_n) + R_{ \geq \frac{dn}{n-1}}  \, .$$
\end{corollary}
\begin{proof} This follows from Theorem \ref{hein} with the help of the
restriction theorems and the numerical formula for tight closure
from the introduction. \end{proof}

\begin{example}
We consider the case of $n$ generic polynomials of degree $d=30$, $3 \leq n \leq 31$.
The following table shows how the degree bound behaves as $n$ varies (we only list $n$ if the
degree bound drops).

\smallskip

\begin{tabular}{c|c|c|c|c|c|c|c|c|c}
$n$ (number of gen. generators) &$3$& $4$ &   $5$ & $6$ & $7$ & $9$ &$11$ & $16$ & $31$
\cr \hline
$\frac{n}{n-1}d$ (degree bound) &
$45$ & $40$ &   $37.5$ & $36$ & $35$ & $33.75$ &$33$  & $32$ & $31$
\end{tabular}

\smallskip

\end{example}

\bibliographystyle{plain}
\bibliography{bibliothek}

\medskip
Holger Brenner, Department of Pure Mathematics, University of
Sheffield, Hicksbuilding, Hounsfield Road, Sheffield S3 7RH, United
Kingdom, e-mail: H.Brenner@sheffield.ac.uk

\markboth{GEORG HEIN}{APPENDIX A: SEMISTABILITY OF THE GENERAL SYZYGY BUNDLE}

\renewcommand  {\ker }  {\operatorname{ker}}

\parindent0cm

\newcommand{\pdop}{{\mathbb P}}
\newcommand{\Ocal}{{\O}}

\renewcommand{\Hom}{{\rm Hom}}
\renewcommand{\proof}{{\bf Proof: }}
\renewcommand{\rk}{{\rm rk}}

\appendix
\section{Semistability of the general Syzygy bundle}

{\hspace{33ex} by Georg Hein}\\

In this appendix we prove three
results about
the (semi)stability of a syzygy bundle. Theorem \ref{app1}
implies Theorem 8.6 of this article.
Here we show stability (resp. semistability) by showing that the restriction of
a sheaf to a given curve is stable (resp. semistable). In theorem \ref{app1} we
use an elliptic curve. This gives us the least restrictive conditions on the integer
parameters $n$ and $d$. However, we can not show stability, because there exist no
stable vector bundles of given rank $r$ and degree $d$ on an elliptic curve unless
$r$ and $d$ are coprime.

Thus, to obtain slope stable coherent sheaves we have to consider curves of
genus greater than 1. This is done in theorems \ref{app2} and \ref{app3}.
It should be remarked that the kernel of a morphisms
$\varphi: \Ocal_{\pdop^N}^{\oplus n} \to \Ocal_{\pdop^N}(d)$
is no vector bundle for $n \leq N$. However, even in these cases we can deduce
(semi)stability.

The strategy of all proofs is as follows:
\begin{itemize}
\item[(1)] We take a suitable (semi)stable sheaf $G$ on a curve $C \subset \pdop^N$.
\item[(2)] We show that there exits a short exact sequence

$$0 \to G \to \Ocal_C^{\oplus n} \stackrel{\bar \varphi}{\rightarrow} \Ocal_{\pdop^N}(d)|_C \to 0 \, .$$

\item[(3)] We show that $\bar \varphi$ is the restriction of a morphisms
$\Ocal_{\pdop^N}^{\oplus n} \stackrel{\varphi}{\rightarrow} \Ocal_{\pdop^N}(d)$ to the curve $C$.
\item[(4)] Now the kernel $F=\ker(\varphi)$ is a coherent sheaf on $\pdop^N$ which is a
vector bundle in an open set containing the curve $C$.
\item[(5)] This implies (see \cite{Sha}) that the restriction of $F$ to the generic curve
in the Hilbert scheme of curves is (semi)stable.
\item[(6)] From that we eventually conclude that $F$ is (semi)stable, because the
restriction of an unstable sheaf to the generic curve in $\pdop^N$ is unstable too.
\vspace{-1ex}\end{itemize}

To show (3) it is sufficient to take projectively normal curves $C \subset \pdop^N$.
We use the theorem of Castelnuovo, Mattuck and Mumford which states that on a curve
$C$ of genus $g_C$ every line bundle $L$ of degree $\deg(L) > 2 \cdot g_C$ is normally
generated (see \cite {GL}). This implies that the embedding $C \to \pdop(H^0(L))$ is
projectively normal.

\begin{theorem}\label{app1}
Let $E \subset \pdop^N$ be a smooth projective elliptic curve embedded by a complete
linear system of degree $N+1$. If the integers $n$ and $d$ satisfy $2 \leq n  \leq d(N+1)$,
then the kernel of a general
morphism $\varphi \in \Hom(\Ocal_{\pdop^N}^{\oplus n}, \Ocal_{\pdop^N}(d) )$ is a semistable
vector bundle when restricted to $E$.
This implies that $\ker(\varphi)$ is a slope semistable coherent sheaf for $\varphi$ generic.
\end{theorem}
\proof
Let $F$ be a semistable vector bundle on the elliptic curve $E$ with $\rk(F)=n-1$ and
$\det(F) \cong \Ocal_{\pdop^N}(d)|_E$. This implies $\deg(F)=d(N+1)$.
The existence of such a vector bundle follows from Atiyah's work \cite{Ati}.
The inequality $n \leq d(N+1)$ implies that $\mu(F) = \frac{\deg(F)}{\rk(F)} > 1$.

Let $P \in E$ be an arbitrary geometric point of $E$.
We consider the  vector bundle $F(-P) = F \otimes \Ocal_E(-P)$.
We compute that the slope $\mu(F(-P))=\mu(F)-1 >0$.
This implies that $H^1(E,F(-P))=0$. Thus, we conclude
from the long exact cohomology sequence associated to
$0 \to F(-P) \to F \to F \otimes k(P) \to 0$
that $F$ is globally generated in the point $P$.
We eventually obtain the surjectivity of the evaluation map
$H^0(E,F) \otimes \Ocal_E \to F$.

By the Riemann-Roch theorem we have $h^0(F) = d(N+1) \geq n$.
Suppose now that  $h^0(F) > n$ holds.
We claim that for a general $n$ dimensional subspace $V \subset H^0(E,F)$
the evaluation morphism ${\rm ev}_V: V\otimes \Ocal_E \to F$ is surjective.
This is done by a dimension count.
The dimension of the Grassmannian variety of all $n$ dimensional subspaces
of $H^0(E,F)$ is $n(h^0(F)-n)$.
Next we consider a surjection $F \stackrel{\alpha}{\rightarrow} k(P)$ and denote its kernel by $F'$.
Since $F$ is globally generated $h^0(F')=h^0(F)-1$.
We deduce that the Grassmannian of all $n$ dimensional subspaces $V$ of $H^0(E,F)$,
such that the image of the evaluation map ${\rm ev}_V$ is contained in $F'$,
is of dimension $n(h^0(F)-n-1)$.
Since the surjections $F \to k(P)$ are parametrized by $\pdop(F)$, and $\dim(\pdop(F))=\rk(F)=n-1$,
we conclude the claim.

Now we take a surjection $\beta:\Ocal_E^{\oplus n} \to F$.
The kernel of this surjection is the line bundle $\Ocal_{\pdop^N}(-d)|_E$.
Thus, considering the dual of $\beta$ we obtain the following short
exact sequence of semistable vector bundles on $E$:
$$0 \lra F^\dual \stackrel{\beta^{\dual}}{\lra} \O_E^{\oplus n}
\stackrel{\bar{\varphi}}{\lra} \O_{\PP^N} (d)|_{E} \lra 0 $$

If we can show that the surjection $\bar \varphi$ is the restriction of a homomorphism
$\varphi \in \Hom(\Ocal_{\pdop^N}^{\oplus n}, \Ocal_{\pdop^N}(d) )$, then our
theorem is proven. In order to conclude our proof, we have to show the surjectivity
of the restriction map
$\Hom(\Ocal_{\pdop^N}^{\oplus n}, \Ocal_{\pdop^N}(d) ) \to
\Hom(\Ocal_{E}^{\oplus n}, \Ocal_{\pdop^N}(d)|_E )$
which is equivalent to the surjectivity of $H^0(\Ocal_{\pdop^N}(d)) \to H^0(\Ocal_{\pdop^N}(d)|_E)$.
However, this is fulfilled since $E$ is projectively normal.
{\hfill $\Box$}

\begin{theorem}\label{app2}
Let $C$ be a smooth quartic in $\pdop^2_k$.
If the integers $n$ and $d$ fulfill the inequality $2 \leq n  \leq \frac{4}{5}d +1$,
then the kernel of a general
morphism $\varphi \in \Hom(\Ocal_{\pdop^2}^{\oplus n}, \Ocal_{\pdop^2}(d) )$ is a stable
vector bundle when restricted to $C$.
This implies that $\ker(\varphi)$ is a slope stable coherent sheaf for a general morphism $\varphi$.
\end{theorem}

\proof
The only new ingredient in our proof is the existence of stable vector bundles
with given determinant on the curve $C$ of genus 3.
This may be deduced from \cite{Ses}.
Indeed, we need a rank $n-1$ stable vector bundle $F$ of determinant $\omega_C^{\otimes d}$.
The slope of $F$ is
$$\mu(F) = \frac{\deg (\omega^{\otimes d} )}{n-1} = \frac{4d}{n-1} \geq 5 \,
.$$
This implies the global
generatedness of $F$ and we can proceed as in the above proof, because $C$ is
projectively normal.
{ \hfill $\Box$}

\begin{theorem}\label{app3}
Let $C \subset \pdop^N$ be a smooth curve of genus two embedded by a complete linear system
of degree $N+2$ for $N \geq 3$. If the integers $n$ and $d$ suffice $2 \leq n \leq \frac{N+2}{3}d+1$,
then the restriction of the kernel of a general morphism $\varphi \in
\Hom(\Ocal_{\pdop^N}^{\oplus n}, \Ocal_{\pdop^N}(d) )$
to $C$ is a stable vector bundle. Thus, $\ker(\varphi)$ is a slope stable coherent sheaf on $\pdop^N$.
\end{theorem}

\proof
As in the proof of theorem \ref{app2}
we have stable vector bundles $F$ with given determinant on $C$.
Since $C$ is of genus two, every stable vector bundle $F$ with $\mu(F) \geq 3$ is globally generated.
The projective normality of $C$ is deduced again by the Castelnuovo, Mattuck, Mumford theorem.
{ \hfill $\Box$}

Georg Hein, Humboldt Universit\"at Berlin, Institut f\"ur Mathematik, Burgstr. 26, 10099 Berlin,
hein@mathematik.hu-berlin.de
\end{document}